# Adaptive estimation of linear functionals in the convolution model and applications

C. BUTUCEA[1] and F. COMTE[2]

[1]*Laboratoire Paul Painlevé – UMR CNRS 8524, Université Lille 1, 59655 Villeneuve d'Ascq Cedex, France. E-mail: Cristina.Butucea@math.univ-lille1.fr*
[2]*MAP5, UMR CNRS 8145, Université Paris Descartes, 45 rue des Saints-Pères, 75270 Paris Cedex 06, France. E-mail: fabienne.comte@parisdescartes.fr*

We consider the model $Z_i = X_i + \varepsilon_i$, for i.i.d. $X_i$'s and $\varepsilon_i$'s and independent sequences $(X_i)_{i\in\mathbb{N}}$ and $(\varepsilon_i)_{i\in\mathbb{N}}$. The density $f_\varepsilon$ of $\varepsilon_1$ is assumed to be known, whereas the one of $X_1$, denoted by $g$, is unknown. Our aim is to estimate linear functionals of $g$, $\langle\psi, g\rangle$ for a known function $\psi$. We propose a general estimator of $\langle\psi, g\rangle$ and study the rate of convergence of its quadratic risk as a function of the smoothness of $g$, $f_\varepsilon$ and $\psi$. Different contexts with dependent data, such as stochastic volatility and AutoRegressive Conditionally Heteroskedastic models, are also considered. An estimator which is adaptive to the smoothness of unknown $g$ is then proposed, following a method studied by Laurent *et al.* (Preprint (2006)) in the Gaussian white noise model. We give upper bounds and asymptotic lower bounds of the quadratic risk of this estimator. The results are applied to adaptive pointwise deconvolution, in which context losses in the adaptive rates are shown to be optimal in the minimax sense. They are also applied in the context of the stochastic volatility model.

## 1. Introduction

We consider the convolution model

$$Z_i = X_i + \varepsilon_i. \tag{1}$$

The sequences $(X_i)_{i\in\mathbb{N}}$ and $(\varepsilon_i)_{i\in\mathbb{N}}$ are independent sequences of real valued random variables. The $X_i$ are i.i.d. with unknown density $g$, the $\varepsilon_i$ are i.i.d. with known density $f_\varepsilon$. The Fourier transform of a function $u \in \mathbb{L}^1(\mathbb{R})$ is denoted by $u^*(x) = \int e^{ixt} u(t)\,dt$. The smoothness of $f_\varepsilon$ is described by parameters $\gamma$, $\alpha$, $\rho$ in the following assumption:

There exist non-negative numbers $\kappa_0$, $\kappa_0'$, $\gamma$, $\alpha$ and $\rho$ such that $f_\varepsilon^*$ satisfies

$$\kappa_0(x^2+1)^{-\gamma/2}\exp\{-\alpha|x|^\rho\} \le |f_\varepsilon^*(x)| \le \kappa_0'(x^2+1)^{-\gamma/2}\exp\{-\alpha|x|^\rho\}, \tag{2}$$

with $\gamma > 1$ when $\rho = 0$. If either $\alpha = 0$ or $\rho = 0$, we set $(\alpha, \rho) = (0, 0)$. Since $f_\varepsilon$ is known, the constants $\alpha, \rho, \kappa_0, \kappa_0'$ and $\gamma$ defined in (2) are known.

When $\rho = 0$ in (2), the errors are called ordinary smooth errors. When $\alpha > 0$ and $\rho > 0$, they are called supersmooth. The standard examples for supersmooth densities are Gaussian or Cauchy distributions (supersmooth of order $\gamma = 0, \rho = 2$ and $\gamma = 0, \rho = 1$, respectively). An example of an ordinary smooth density is the Laplace distribution ($\rho = 0 = \alpha$ and $\gamma = 2$).

In this context, many papers have studied the deconvolution problem. Many different strategies have been developed in order to estimate the distribution $g$ of the unobserved $X_i$, when $g$ is assumed to belong to some smoothness class defined by

$$\mathcal{S}(b,a,r,L) = \left\{ g \text{ such that } \int_{-\infty}^{+\infty} |g^*(x)|^2 (x^2+1)^b \exp\{2a|x|^r\}\,\mathrm{d}x \le 2\pi L \right\}, \tag{3}$$

where $b, a, r, L$ are some unknown non-negative numbers, such that $b > 1/2$ when $r = 0$. If either $a = 0$ or $r = 0$, we set $(a, r) = (0, 0)$ and we say that the density is ordinary smooth. When both $a$, $r > 0$, we call the density supersmooth.

In this paper, we are interested in the problem of estimating $\theta(g) = \langle \psi, g\rangle = \mathbb{E}(\psi(X_1))$ in model (1), where $\psi$ is a known integrable function with respect to the probability measure associated to $g$. To study the rates of convergence of our estimators, we have to take into account the smoothness of the function $\psi$. Thus $\psi$ is assumed to satisfy:

$$\forall x \in \mathbb{R} \qquad |\psi^*(x)|^2 \le C_\psi (x^2+1)^{-B} \exp(-2A|x|^R). \tag{4}$$

The parameters $A$ and $R$ are non-negative real numbers, and $B$ is non-negative or such that $\psi^* g^*$ is integrable. In particular, they can be zero if $g^*$ is integrable. We work under the convention that if either $A = 0$ or $R = 0$, then we set $(A, R) = (0, 0)$.

We exhibit the whole range of the rates of convergence for estimators of the functional $\theta(g)$, depending on the parameters in (2)–(4). To the best of our knowledge, this general rate description is new. We also extend the result to different dependency contexts, in view of applications to particular hidden Markov models or AutoRegressive Conditionally Heteroskedastic-type models.

The upper bounds for the rates follow from a squared-bias/variance compromise. To obtain this compromise, we have to choose a smoothing parameter which depends on unknown quantities. Therefore, a data driven model selection type procedure is proposed. It is based on minimization of a penalized estimated criterion, which is different from the one intensively studied for mean integrated squared errors. The difficulty here lies in finding an adequate criterion for the setting of a linear functional and mean squared error. The proposed procedure is inspired by Laurent *et al.* [24]. We give upper bounds for this adaptive method, with particular interest in the cases where a loss in the rate appeared with respect to the non-adaptive estimator.

In the particular case of pointwise estimation, adaptive estimation in the direct problem (i.e., when the $X_i$ are observed without noise) has been widely studied in the context of the Gaussian white noise and regression models, see, for example, Lepski [26], Tsybakov [28], Cai and Low [7, 8] (for more general linear functionals), Artiles and Levit [2], Laurent *et al.* [24] and, in the context of density models, Lepski and Levit [25], Butucea [3] and Artiles [1]. For the model of Gaussian sequences Golubev and Levit [21] and Golubev [20] considered adaptive estimation of linear functionals in both direct and inverse problems. In the Gaussian white noise model Goldenshluger [18] and Goldenshluger and Pereverzev [19] considered pointwise estimation for the inverse problem on classes of functions similar to $S(b,0,0,L)$. Their adaptive procedure is based on Lepski's procedure. Note also that in some particular inverse problems the pointwise adaptive estimation was solved by Klemelä and Tsybakov [22] for the Riesz transform and by Cavalier [10] for the tomography problem. To the best of our knowledge, we present the first work on adaptive estimation of general functionals of the form $\int \psi g$ in the context of indirect observation (1).

We do not study optimality in the very general case: this would be very technical. But we study as a first application the particular case of pointwise density deconvolution. This case corresponds to $\psi^*(t) = \mathrm{e}^{\mathrm{i}tx_0}$, which satisfies (4), meaning that we can choose $\psi$ as the Dirac measure at $x_0$. This makes sense in our problem because the definition of our estimator involves only $\psi^*$. We recover in this particular case the upper bound rates obtained by Fan [16], Cator [9], Butucea [4] and Butucea and Tsybakov [6]. Moreover, we prove the optimality in the minimax sense of the loss due to adaptation for Sobolev smooth and supersmooth densities in the presence of ordinary smooth noise and for supersmooth densities in the presence of supersmooth noise with $r \geq \rho$ and $0 < \rho \leq 1$ (in the case $r < \rho$ no loss occurs, while the case $r \geq \rho$ and $1 < \rho < 2$ is still open). As a by-product we also prove in the last case that the rate of our estimator is optimal in the minimax sense, which was not yet known in the literature.

Our estimation method is also illustrated for the discrete stochastic volatility model, where derivatives of the Laplace transform of the volatility are estimated with good rates.

The plan of the paper is as follows: In Section 2 we define the estimators and we compute upper bounds for their mean squared error. In Section 3 the adaptive procedure is detailed. Both independent and $\beta$-mixing contexts are studied. In Section 4, two applications of our general results are detailed. Section 4.1 shows the application of the results to adaptive pointwise deconvolution, upper bounds are deduced from Section 3 and the associated lower bounds are proven when a loss occurs. Section 4.2 presents an application to the context of the stochastic volatility model. Most proofs are gathered in Section 5.

## 2. Risk bound for the estimator

We denote by $\langle \cdot,\cdot \rangle$ the $\mathbb{L}^2$-scalar product ($\langle u,v\rangle = \int u(x)\bar{v}(x)\,\mathrm{d}x$), by $\star$ the convolution product of functions ($u \star v(x) = \int u(t)v(t-x)\,\mathrm{d}t$) and by $u^*$ the Fourier transform of $u \in \mathbb{L}^1(\mathbb{R})$: $u^*(x) = \int \mathrm{e}^{\mathrm{i}tx}u(t)\,\mathrm{d}t$.

Recall that we want to estimate $\theta(g) = \langle \psi, g \rangle = \mathbb{E}(\psi(X_1))$ where $X_1$ follows model (1) and is unobserved. Only the $Z_i$, for $i = 1, \ldots, n$ are available. In what follows we assume that

$$\begin{cases} \psi g,\ \psi \text{ and } \psi^* g^* \text{ belong to } \mathbb{L}_1(\mathbb{R}), \\ f_\varepsilon \text{ belongs to } \mathbb{L}_2(\mathbb{R}) \text{ and is such that } \forall x \in \mathbb{R}\ f_\varepsilon^*(x) \neq 0. \end{cases} \tag{5}$$

Note that the square integrability of $f_\varepsilon$ requires that $\gamma > 1/2$ when $\rho = 0$ in (2).

Moreover, we generalize these results to distributions having Fourier transform such that $\int \psi^* g^* < \infty$. For example, we estimate $g(x_0)$ for some fixed $x_0$ when we take $\psi = \delta_{x_0}$, the Dirac measure at $x_0$, having Fourier transform equal to $\psi^*(t) = \mathrm{e}^{-\mathrm{i}tx_0}$. For estimating the derivatives $g^{(k)}(x_0)$, when they exist, we consider $\psi$ such that $\psi^*(t) = (-\mathrm{i}t)^k \mathrm{e}^{-\mathrm{i}tx_0}$.

## 2.1. The estimator

We write, using (5), $\langle \psi, g \rangle = (1/2\pi)\langle \psi^*, g^* \rangle = (1/2\pi)\langle \psi^*, f_Z^*/f_\varepsilon^* \rangle$. Replacing $f_Z^*(t)$ by its empirical version $(1/n)\sum_{k=1}^n \mathrm{e}^{\mathrm{i}tZ_k}$ leads to the estimator

$$\hat{\theta} = \frac{1}{2\pi n} \sum_{k=1}^n \int \mathrm{e}^{\mathrm{i}tZ_k} \frac{\psi^*(-t)}{f_\varepsilon^*(t)}\, \mathrm{d}t. \tag{6}$$

This estimator is explicit and seems attractive. Unfortunately, the integral diverges for many choices of $f_\varepsilon^*$; for instance, $\varepsilon$ is a Gaussian noise. To overcome such issues, we suggest regularization and take the following estimator of $\theta(g)$:

$$\hat{\theta}_m = \frac{1}{2\pi n} \sum_{k=1}^n \int_{|t| \leq \pi m} \mathrm{e}^{\mathrm{i}tZ_k} \frac{\psi^*(-t)}{f_\varepsilon^*(t)}\, \mathrm{d}t, \tag{7}$$

where $m$ is an integer.

***Remark 2.1.*** Let $\hat{g}_m$ denote the projection estimator of $g$ defined in Comte *et al.* [14]. Then we can prove that $\hat{\theta}_m = \langle \hat{g}_m, \psi \rangle$.

## 2.2. Risk bounds and rates for i.i.d. variables $X_i$'s

If (5) holds and if, moreover, $\psi^*(-\cdot)/f_\varepsilon^*$ is integrable, then $m = +\infty$ can be chosen and the estimator $\hat{\theta} = \hat{\theta}_\infty$ is unbiased and has a parametric rate.

Otherwise, we have $\mathbb{E}(\theta - \hat{\theta}_m)^2 = b^2(\hat{\theta}_m) + \mathrm{Var}(\hat{\theta}_m)$ with $b(\hat{\theta}_m) = \theta - \mathbb{E}(\hat{\theta}_m)$. As $\mathbb{E}(\hat{\theta}_m) = (2\pi)^{-1} \int_{|t| \leq \pi m} g^*(t)\psi^*(-t)\, \mathrm{d}t$, we obtain

$$\begin{aligned} b(\hat{\theta}_m) &= \frac{1}{2\pi}\left( \int g^*(t)\psi^*(-t)\, \mathrm{d}t - \int_{|t| \leq \pi m} g^*(t)\psi^*(-t)\, \mathrm{d}t \right) \\ &= \frac{1}{2\pi} \int_{|t| \geq \pi m} g^*(t)\psi^*(-t)\, \mathrm{d}t. \end{aligned} \tag{8}$$

Under (5), $b(\hat\theta_m)$ tends to 0 when $m$ grows to infinity.

For the variance term, write

$$\operatorname{Var}(\hat\theta_m) = \frac{1}{4\pi^2 n}\operatorname{Var}\left(\int_{-\pi m}^{\pi m} \mathrm{e}^{\mathrm{i}uZ_1}\frac{\psi^*(-u)}{f_\varepsilon^*(u)}\,\mathrm{d}u\right).$$

First, the following bound holds:

$$\operatorname{Var}(\hat\theta_m) \le \frac{1}{4\pi^2 n}\mathbb{E}\left(\left|\int_{-\pi m}^{\pi m} \mathrm{e}^{\mathrm{i}uZ_1}\frac{\psi^*(-u)}{f_\varepsilon^*(u)}\,\mathrm{d}u\right|^2\right) \le \frac{1}{4\pi^2 n}\left(\int_{-\pi m}^{\pi m}\frac{|\psi^*(-u)|}{|f_\varepsilon^*(u)|}\,\mathrm{d}u\right)^2.$$

Next, the variance can also be bounded as follows:

$$\begin{aligned}\operatorname{Var}(\hat\theta_m) &= \frac{1}{4\pi^2 n}\int_{-\pi m}^{\pi m}\int_{-\pi m}^{\pi m}(f_Z^*(u-v) - f_Z^*(u)f_Z^*(-v))\frac{\psi^*(-u)\psi^*(v)}{f_\varepsilon^*(u)f_\varepsilon^*(-v)}\,\mathrm{d}u\,\mathrm{d}v\\ &\le \frac{1}{4\pi^2 n}\int_{-\pi m}^{\pi m}\int_{-\pi m}^{\pi m}\frac{\psi^*(-u)\psi^*(-v)}{f_\varepsilon^*(u)f_\varepsilon^*(-v)}f_Z^*(u-v)\,\mathrm{d}u\,\mathrm{d}v.\end{aligned}$$

We use the Cauchy–Schwarz inequality and Fubini's theorem:

$$\begin{aligned}\operatorname{Var}(\hat\theta_m) &\le \frac{1}{4\pi^2 n}\left(\int_{-\pi m}^{\pi m}\int_{-\pi m}^{\pi m}\left|\frac{\psi^*(-u)}{f_\varepsilon^*(u)}\right|^2|f_Z^*(u-v)|\,\mathrm{d}u\,\mathrm{d}v\right.\\ &\qquad\left.\times\int_{-\pi m}^{\pi m}\int_{-\pi m}^{\pi m}\left|\frac{\psi^*(v)}{f_\varepsilon^*(-v)}\right|^2|f_Z^*(u-v)|\,\mathrm{d}u\,\mathrm{d}v\right)^{1/2}\\ &\le \frac{1}{4\pi^2 n}\int_{|u|\le\pi m}\frac{|\psi^*(-u)|^2}{|f_\varepsilon^*(u)|^2}\,\mathrm{d}u\int|f_Z^*(x)|\,\mathrm{d}x.\end{aligned}$$

Note that since $\psi$ is a real valued function we have $|\psi^*(-t)| = |\psi^*(t)|$. As $\int|f_Z^*(x)|\,\mathrm{d}x \le \int|f_\varepsilon^*(x)|\,\mathrm{d}x$, we have the following result:

**Proposition 2.1.** *Assume that $C_\varepsilon = \int|f_\varepsilon^*(x)|\,\mathrm{d}x < +\infty$, and let $\hat\theta_m$ be defined by (7). Then, under (5),*

$$\mathbb{E}(\theta-\hat\theta_m)^2 \le \left(\frac{1}{2\pi}\int_{|t|\ge\pi m}|g^*(t)\psi^*(t)|\,\mathrm{d}t\right)^2 + \frac{1}{4\pi^2 n}\min\left\{C_\varepsilon\int_{-\pi m}^{\pi m}\frac{|\psi^*|^2}{|f_\varepsilon^*|^2}, \left(\int_{-\pi m}^{\pi m}\frac{|\psi^*|}{|f_\varepsilon^*|}\right)^2\right\}.$$

Note that we also have $\int|f_Z^*(x)|\,\mathrm{d}x \le \|f_\varepsilon^*\|\|g^*\| = 2\pi\|f_\varepsilon\|\|g\|$, if $f_\varepsilon$ and $g$ are square integrable.

***Remark 2.2.*** If, in addition,

$$\int|\psi^*(x)/f_\varepsilon^*(x)|^2\mathrm{d}x < +\infty, \tag{9}$$

then the variance of $\hat{\theta}_m$ is of order $1/n$ and the estimator can reach the parametric rate, for $m$ large enough. Note that a condition like $\int |\psi^*(x)/f_\varepsilon^*(x)|\,\mathrm{d}x < \infty$ (which ensures that (6) is well defined) is generally stronger than (9) as convergence problems lie only near infinity. Moreover, such conditions are fulfilled if $\psi^*$ decreases faster than $f_\varepsilon^*$ near infinity, which corresponds to the intuitive idea that $\psi$ is a smoother function than $f_\varepsilon$. For example, this happens if $\psi$ is supersmooth and $f_\varepsilon$ is ordinary smooth.

Thus we can study the rates that can be deduced from the upper bounds of Proposition 2.1, as a function of the smoothness parameters of the three functions involved, $g$, $\psi$, $f_\varepsilon$. To do so, let us assume that $\psi$ satisfies (4), that $g$ belongs to $\mathcal{S}(b,a,r,L)$ as defined by (3) and that $f_\varepsilon^*$ fulfills (2). Then, use (8), (3) and (4) to get

$$\begin{aligned} b^2(\hat{\theta}_m) &\le \left| \int_{|x|\ge \pi m} |g^*(x)|(1+x^2)^{b/2}\exp(a|x|^r)(|\psi^*(x)|(1+x^2)^{-b/2}\exp(-a|x|^r))\,\mathrm{d}x \right|^2 \\ &\le \int_{|x|\ge \pi m} |g^*(x)|^2(1+x^2)^b \exp(2a|x|^r)\,\mathrm{d}x \\ &\quad \times \int_{|x|\ge \pi m} |\psi^*(x)|^2(1+x^2)^{-b}\exp(-2a|x|^r)\,\mathrm{d}x \\ &\le LC\int_{|x|\ge \pi m} (1+x^2)^{-b-B}\exp(-2a|x|^r - 2A|x|^R)\,\mathrm{d}x \\ &\le C_1 m^{-2b-2B-\max(r,R)+1}\exp(-2a(\pi m)^r - 2A(\pi m)^R). \end{aligned}$$

On the other hand, the noise plays an important role on the variance of the estimator:

- Case (I): If ($\rho = R = 0, \gamma < B - 1/2$) or ($\rho = R > 0, \alpha = A, \gamma < B - 1/2$) or ($\rho = R$, $\alpha < A$) or ($\rho < R$), then $\mathrm{Var}(\hat{\theta}_m) \le C' n^{-1}$;
- Case (II): If ($\rho = R = 0, \gamma = B - 1/2$) or ($\rho = R > 0, \alpha = A, \gamma = B - 1/2$), then $\mathrm{Var}(\hat{\theta}_m) \le C' \ln(m) n^{-1}$;
- Case (III): If ($\rho = R = 0, \gamma > B - 1/2$) or ($\rho = R > 0, \alpha = A, \gamma > B - 1/2$), then $\mathrm{Var}(\hat{\theta}_m) \le C' m^{2\gamma - 2B + 1} n^{-1}$;
- Case (IV): If ($\rho > R$) or ($\rho = R > 0, \alpha > A$), then

$$\mathrm{Var}(\hat{\theta}_m) \le C' n^{-1} m^{2\gamma - 2B + 1 - \rho + (1-\rho)_+} \mathrm{e}^{2\alpha(\pi m)^\rho - 2A(\pi m)^R}.$$

We summarize in Table 1 the scenarios that arise when one minimizes over $m$ the sum of the upper bounds on the bias and the variance. Let $a \vee b = \max\{a, b\}$. Note that in cases (8) and (9) the rate is given by

$$\begin{aligned} v_n = \min_m \Big\{ & C_B m^{-2b-2B+1-r\vee R}\mathrm{e}^{-2a(\pi m)^r - 2A(\pi m)^R} \\ & + m^{2\gamma-2B+1-\rho+(1-\rho)_+}\mathrm{e}^{2\alpha(\pi m)^\rho - 2A(\pi m)^R}\frac{1}{n}\Big\}. \end{aligned} \tag{10}$$

**Table 1.** Upper bounds for the minimax rates of convergence, $\delta_1 = (2\gamma - 2B + 1)/\{r \vee R\}$, $\delta_2 = (b + B - 1/2)/(b + \gamma)$ and $\delta_3 = (2(b + B) - 1)/\rho$

| Parameters | | | Rates | Adaptive rates |
|---|---|---|---|---|
| $\rho < R$ | | | (1) $n^{-1}$ | |
| $\rho = R$, $\alpha < A$ | | | (2) $n^{-1}$ | |
| $(\rho = R = 0)$ or $(\rho = R > 0,$ $\alpha = A)$ | $\gamma < B - 1/2$ | | (3) $n^{-1}$ | |
| | $\gamma = B - 1/2$ | $r \vee R > 0$ | (4) $(\ln\ln n)n^{-1}$ | $(\ln\ln n)^2 n^{-1}$ |
| | | $r \vee R = 0$ | (5) $(\ln n)n^{-1}$ | $(\ln n)^2 n^{-1}$ |
| | $\gamma > B - 1/2$ | $r \vee R > 0$ | (6) $(\ln n)^{\delta_1} n^{-1}$ | $(\ln\ln(n)\ln n)^{\delta_1} n^{-1}$ |
| | | $r \vee R = 0$ | (7) $n^{-\delta_2}$ | $(n/\ln(n))^{-\delta_2}$ |
| $(\rho = R > 0,$ $\alpha > A)$ | | | (8) $v_n$ in (10) | $v_n(\ln(n))^{\delta_4}, 0 \le \delta_4 \le 1$ |
| $\rho > R$ | | $r \vee R > 0$ | (9) $v_n$ in (10) | $v_n(\ln(n))^{\delta_4}, 0 \le \delta_4 \le 1$ |
| | | $r \vee R = 0$ | (10) $\ln(n)^{-\delta_3}$, | $\ln(n)^{-\delta_3}$ |

These rates are strictly faster than $(\ln(n))^{-\lambda_1}$, that is, $v_n = o((\ln(n))^{-\lambda_1})$ for any $\lambda_1 > 0$, and generally slower than $n^{-\lambda_2}$, $\lambda_2 > 0$ (negative powers of $n$ can be obtained). For precise (but cumbersome) formulae in similar cases, we refer to Lacour [23]. We give in Section 5 the orders of the $m$ associated to the rates.

***Remark 2.3.*** Different optimal choices of $m$ depend on the unknown parameters related to $g$ (see Section 5.1), hence the interest in an automatic selection procedure for $m$.

## 2.3. Extension to mixing contexts

In view of applications, it is natural to study the robustness of our method when the variables $X_i$ are $\beta$-mixing. To be more precise, two dependence contexts are considered.

- (D1) In Model (1), the sequences $(X_i)$ and $(\varepsilon_i)$ are independent and the $\varepsilon_i$ are i.i.d. The sequence $(X_i)$ is strongly stationary and $\beta$-mixing, with $\beta$-mixing coefficients denoted by $(\beta_k)_k$.
- (D2) In Model (1), the $\varepsilon_i$ are i.i.d. and, for any given $i$, $X_i$ and $\varepsilon_i$ are independent (but the sequences $(X_i)$ and $(\varepsilon_i)$ are not independent). The sequence $(Z_i, X_i)_{i\in\mathbb{Z}}$ is strongly stationary and $\beta$-mixing, with $\beta$-mixing coefficients denoted by $(\beta_k)_k$.

Context (D1) encompasses the case of particular hidden Markov models, when the noise is additive and $(X_i)$ is a $\beta$-mixing Markov process. As many Markov chain models or other standard models can be proved to have such mixing properties (see Doukhan [15] for a large number of examples and study of their mixing properties), this means that our results can be applied to many classical models. In that case, we can prove the following result:

**Proposition 2.2.** *Consider the model (1) under* (D1) *with moreover* $\sum_{k\geq 0}\beta_k < +\infty$. *Assume that* $C_\varepsilon = \int |f_\varepsilon^*(x)|\,\mathrm{d}x < +\infty$. *Let* $\hat\theta_m$ *be defined by (7). Then*

$$\mathbb{E}(\theta-\hat\theta_m)^2 \leq \left(\frac{1}{2\pi}\int_{|t|\geq \pi m}|g^*(t)\psi^*(t)|\,\mathrm{d}t\right)^2 + \frac{C_\varepsilon}{4\pi^2 n}\min\left\{\int_{-\pi m}^{\pi m}\frac{|\psi^*|^2}{|f_\varepsilon^*|^2}, \left(\int_{-\pi m}^{\pi m}\frac{|\psi^*|}{|f_\varepsilon^*|}\right)^2\right\} + \frac{2(\int_{|t|\leq \pi m}|\psi^*|(t)\,\mathrm{d}t)^2\sum_{k\geq 0}\beta_k}{n}. \tag{11}$$

In particular, if $K_\psi := \int |\psi^*(t)|\,\mathrm{d}t < +\infty$, then the last term in the right-hand side of (11) is of order $\mathrm{O}(1/n)$. Moreover, in any case, we have in (11),

$$\left(\int_{|t|\leq \pi m}|\psi^*(t)|\,\mathrm{d}t\right)^2 \leq \min\left\{2\pi\|f_\varepsilon\|^2\int_{-\pi m}^{\pi m}\frac{|\psi^*|^2}{|f_\varepsilon^*|^2}, \left(\int_{-\pi m}^{\pi m}\frac{|\psi^*|}{|f_\varepsilon^*|}\right)^2\right\},$$

so that the last term is always less than or equal to the variance term. It follows that the rates, in the context of mixing $X_k$ described by assumption (D1), remain the same as in the independent setting.

Context (D2) is related to ARCH models. Indeed, general ARCH models can be formulated as follows: Let $(\eta_i)$ be an i.i.d. noise sequence.

$$Y_i = \sigma_i\eta_i \qquad \text{with } \sigma_i = F(\eta_{i-1}, \eta_{i-2}, \ldots), \tag{12}$$

for some measurable functions $F$, or

$$Y_i = \sigma_i\eta_i \qquad \text{with } \sigma_i = F(\sigma_{i-1}, \eta_{i-1}) \text{ and } \sigma_0 \text{ independent of } (\eta_i)_{i\geq 0}. \tag{13}$$

Many examples can be found in the literature, and conditions can be given under which the process $(Y_i,\sigma_i)_{i\in\mathbb{Z}}$ is geometrically $\beta$-mixing; we refer to Comte *et al.* [12] for a review of the examples and the references therein. Clearly then, $Z_i = \ln(Y_i^2)$, $X_i = \ln(\sigma_i^2)$ and $\varepsilon_i = \ln(\eta_i^2)$ follow model (1) and satisfy conditions given by (D2). We can prove the following result in this context:

**Proposition 2.3.** *Consider the model (1) under* (D2) *with moreover* $\sum_{k\geq 0}\beta_k < +\infty$. *Assume that* $C_\varepsilon = \int |f_\varepsilon^*(x)|\,\mathrm{d}x < +\infty$. *Let* $\hat\theta_m$ *be defined by (7). Then*

$$\mathbb{E}(\theta-\hat\theta_m)^2 \leq \left(\frac{1}{2\pi}\int_{|t|\geq \pi m}|g^*(t)\psi^*(t)|\,\mathrm{d}t\right)^2$$

$$+\frac{C_\varepsilon}{4\pi^2 n}\min\left\{\int_{-\pi m}^{\pi m}\frac{|\psi^*|^2}{|f_\varepsilon^*|^2},\left(\int_{-\pi m}^{\pi m}\frac{|\psi^*|}{|f_\varepsilon^*|}\right)^2\right\} \tag{14}$$
$$+\frac{2\sum_{k\geq 0}\beta_k}{n}\left(\int_{-\pi m}^{\pi m}|\psi^*|\right)\left(\int_{-\pi m}^{\pi m}\frac{|\psi^*|}{|f_\varepsilon^*|}\right).$$

Thus, the procedure attains the rates of the independent case as soon as, for some constant $C$,

$$\left(\int_{-\pi m}^{\pi m}|\psi^*|\right)\left(\int_{-\pi m}^{\pi m}\frac{|\psi^*|}{|f_\varepsilon^*|}\right)\leq C\min\left\{\int_{-\pi m}^{\pi m}\frac{|\psi^*|^2}{|f_\varepsilon^*|^2},\left(\int_{-\pi m}^{\pi m}\frac{|\psi^*|}{|f_\varepsilon^*|}\right)^2\right\}.$$

This does not hold in general, but in particular cases. For instance, if $f_\varepsilon$ satisfies (2) and if $\psi$ satisfies (4) together with

$$|\psi^*(x)|^2\geq C'_\psi(x^2+1)^{-B}\exp(-2A|x|^R), \tag{15}$$

with either $\gamma>\max(B,1)$ or $(A>0,\rho>0)$, then, under the assumptions of Proposition 2.3,

$$\mathbb{E}(\theta-\hat\theta_m)^2\leq\left(\frac{1}{\pi}\int_{\pi m}^{+\infty}|g^*\psi^*|\right)^2+\frac{K}{4\pi^2 n}\min\left\{\int_{-\pi m}^{\pi m}\frac{|\psi^*|^2}{|f_\varepsilon^*|^2},\left(\int_{-\pi m}^{\pi m}\frac{|\psi^*|}{|f_\varepsilon^*|}\right)^2\right\}, \tag{16}$$

where $K$ is a constant.

It follows from (16) that the rates given in Table 1 are preserved in this $\beta$-mixing context whenever the $\varepsilon_i$ are supersmooth.

Taking $\psi^*(t)=\mathrm{e}^{\mathrm{i}tx_0}$ for any $x_0$ (as in Section 4.1 below) allows one to provide a pointwise density estimator that retains the rate of the independent case if $\gamma>1$. We recover the results obtained by the kernel estimator of van Es *et al.* [29]. Our results are more general since van Es *et al.* [29] only consider a multiplicative Gaussian noise (implying supersmooth $\varepsilon_i$, see Section 4.2) and do not study adaptation (which is not useful in their particular case). Other functionals $\langle\psi,g\rangle$ may be estimated with our procedure.

## 3. Adaptive estimation

Now, we provide a strategy leading to an automatic choice of $m$. Note that such model selection has an interest only in the case $\int|\psi^*/f_\varepsilon^*|=+\infty$ and $\int|\psi^*/f_\varepsilon^*|^2=+\infty$ since otherwise the variance is of order $1/n$ and the rate is parametric. As $\psi$ and $f_\varepsilon$ are assumed to be known, these conditions can be explicitly checked.

Let us describe briefly the heuristics that follow Laurent *et al.* [24]. Let $\theta_m=\mathbb{E}(\hat\theta_m)=(2\pi)^{-1}\int_{-\pi m}^{\pi m}g^*(t)\psi^*(-t)\,\mathrm{d}t$. The approximation of the bias of $(\theta(g)-\theta_m)^2$ is obtained by replacing it by $(\theta_j-\theta_m)^2$ for $j\geq m$, $j$ great enough, and then by $(\hat\theta_j-\hat\theta_m)^2$. This approximation in turn introduces a bias which must be corrected (see $H(j,m)$ below).

The variance term is replaced by a penalty function $\mathrm{pen}(\cdot)$ from $\mathbb{N}$ into $\mathbb{R}^+$. This gives the theoretical criterion

$$Crit(m) = \sup_{j\geq m}(\theta_j - \theta_m)^2 + \mathrm{pen}(m),$$

where $\mathrm{pen}(m)$ has the order of the variance term (see Section 2.2) and its empirical version is

$$\widehat{Crit}(m) = \sup_{j\geq m, j\in\mathcal{M}}[(\hat{\theta}_m - \hat{\theta}_j)^2 - H(j,m)] + \mathrm{pen}(m),$$

where $H(j,m)$ is an additional bias correction and $\mathcal{M}$ is a subset of $\mathbb{N}$. Then, we can define

$$\hat{m} = \inf\left\{m\in\mathcal{M}, \widehat{Crit}(m) \leq \inf_{j\in\mathcal{M}} \widehat{Crit}(j) + \frac{1}{n}\right\} \tag{17}$$

as the model selection procedure. It remains to find $\mathrm{pen}(\cdot)$ and $H(j,m)$ that make the procedure work and give good rates for $\hat{\theta}_{\hat{m}}$.

Recall that $C_\varepsilon = \int |f_\varepsilon^*(x)|\,\mathrm{d}x$. Let $x_m$, be some positive weights to be chosen, and let $\mathrm{a} > 0$. We define:

$$\mathrm{pen}(m) = 4\left(1+\frac{1}{\mathrm{a}}\right)(x_m\sigma_m^2 + x_m^2 c_m^2), \tag{18}$$

where $\sigma_m^2 = \sigma_{0,m}^2$, $c_m = c_{0,m}$, with $\sigma_{j,m}^2$ and $c_{j,m}$ defined by

$$\sigma_{j,m}^2 = \frac{1}{2\pi n}\min\left\{C_\varepsilon \int_{\pi(j\wedge m)\leq|x|\leq\pi(j\vee m)} \left|\frac{\psi^*(x)}{f_\varepsilon^*(x)}\right|^2 \mathrm{d}x,\right.$$
$$\left.\left(\int_{\pi(j\wedge m)\leq|x|\leq\pi(j\vee m)} \frac{|\psi^*(x)|}{|f_\varepsilon^*(x)|}\,\mathrm{d}x\right)^2\right\}$$

and

$$c_{j,m} = \frac{1}{2\pi n}\int_{\pi(j\wedge m)\leq|x|\leq\pi(j\vee m)} \left|\frac{\psi^*(x)}{f_\varepsilon^*(x)}\right|\mathrm{d}x.$$

Let also

$$H(j,m) = 4\left(1+\frac{1}{\mathrm{a}}\right)(x_j\sigma_{j,m}^2 + x_j^2 c_{j,m}^2). \tag{19}$$

We shall prove the following theorem:

**Theorem 3.1.** *Consider model (1) where $(X_i)_{1\leq i\leq n}$ and $(\varepsilon_i)_{1\leq i\leq n}$ are independent sequences of i.i.d. random variables and assume that (5) is fulfilled. Let $\hat{\theta}_{\hat{m}}$ be defined by (7) and (17)–(19) when $\int|\psi^*/f_\varepsilon^*| = +\infty$ and $\int|\psi^*/f_\varepsilon^*|^2 = +\infty$. Then there exists some* *positive constant* $C(\mathrm{a})$ *depending only on some* $\mathrm{a} > 0$*, such that*

$$\mathbb{E}[(\hat{\theta}_{\hat{m}} - \theta)^2] \leq C(\mathrm{a}) \inf_{m \in \mathcal{M}} \left\{ \left( \int_{|x| \geq \pi m} |\psi^*(x) g^*(x)| \, \mathrm{d}x \right)^2 + \mathrm{pen}(m) \right\}$$
$$+ C(\mathrm{a}) \sum_{m \in \mathcal{M}} \mathrm{e}^{-x_m} \omega_m^2 + \frac{C'}{n},$$

*where* $\omega_m^2 = \sigma_m^2 \vee c_m + 2(\sigma_m^2 \vee c_m)^2$ *and* $C'$ *is a constant.*

Theorem 3.1 states that $\hat{\theta}_{\hat{m}}$ leads to an automatic tradeoff between the squared bias term $(\int_{|x| \geq \pi m} |\psi^*(-x) g^*(x)| \, \mathrm{d}x)^2$ and $\mathrm{pen}(m)$, if $x_m$ are chosen so that $\sum_m \mathrm{e}^{-x_m} \omega_m^2 = \mathrm{O}(1/n)$. However, as the main term in $\mathrm{pen}(m)$ is clearly $x_m \sigma_m^2$, where $\sigma_m^2$ is the variance of $\hat{\theta}_m$, $x_m$ represents a loss in the variance (not necessarily in the rate).

Now, let us discuss the possible choices for $x_m$ in order to see what loss occurs, if any, when using the adaptive procedure. We assume that $b + B > 1$, so that we can take $\mathcal{M} = \{1, 2, \ldots, [\sqrt{n}]\}$, where $[\sqrt{n}]$ is the greatest integer less than $\sqrt{n}$. The possible choices for $x_m$ are discussed with respect to the upper bounds on the variance given in Section 2.2:

- Case (II): We take $x_m = 2\ln(m)$ and the rate becomes of order $(\ln\ln(n))^2/n$ instead of $\ln\ln(n)/n$ or of order $\ln^2(n)/n$ instead of $\ln(n)/n$.
- Case (III): We take $x_m = (2\gamma - 2B + 3)\ln(n)$ and the rate becomes of order $\ln\ln(n)\ln^\delta(n)/n$ instead of $\ln^\delta(n)/n$ and of order $(n/\ln(n))^{-[(b+B)-1/2]/(b+\gamma)}$ instead of $n^{-[(b+B)-1/2]/(b+\gamma)}$.
- Case (IV): We take $x_m = 4\alpha(\pi m)^\rho$. There is no loss in case (10) if $\rho > 0$, $r = R = 0$. In the two other cases, (8) and (9), a loss in the variance occurs. If the bias is dominating (if $r > \rho$), there is no loss in the rate. Otherwise, as the optimal $m$ is less than $(\ln(n)/C)^{1/\rho}$, for some $C > 0$, the loss in the rate is at most of order $\mathrm{O}(\ln(n))$. Note that the rate being faster than logarithmic in this case, the loss remains negligible with respect to the rate.

The adaptive rates are given in the last column of Table 1. Let us emphasize that the rates presented in both the second and third columns of Table 1 are new in such a general setup.

Moreover, if we want to extend the adaptive result to the mixing case, we can use the Bernstein inequality given in Doukhan [15] or in Butucea and Neumann [5] provided that the mixing is geometrical. We can prove the following corollary of Theorem 3.1:

**Corollary 3.1.** *Consider model (1) under* (D1) *or under* (D2) *with* $f_\varepsilon$ *satisfying (2) and* $\psi$ *satisfying (4) and (15) with either* $\gamma > \max(B, 1)$ *or* $A, \rho > 0$*, and assume in both cases that* $\beta_k \leq \mathrm{e}^{-ck}$ *for any* $k \in \mathbb{N}$*. Then if (5) is fulfilled and if* $\int |\psi^*(t)| \, \mathrm{d}t < +\infty$, $\int |\psi^*/f_\varepsilon^*| = +\infty$ *and* $\int |\psi^*/f_\varepsilon^*|^2 = +\infty$*, the result of Theorem 3.1 for* $\hat{\theta}_{\hat{m}}$ *defined in the same way holds with* $c_m$, $c_{j,m}$ *replaced by* $2c_m \ln(n)/c$, $2c_{j,m}\ln(n)/c$ *and* $\sigma_m^2$, $\sigma_{j,m}^2$ *multiplied by 2.*

Clearly, the constant $c$ appearing in the $c_m$, $c_{m,j}$ is unknown, but these terms have in general negligible orders when compared to the $\sigma_m^2$, $\sigma_{j,m}^2$. In that case, these terms can be omitted in the definition of the estimator and the procedure does not depend on the mixing coefficients (see the example in Section 4.2).

# 4. Applications

## 4.1. Pointwise estimation

Pointwise estimation of $g$, also called pointwise deconvolution, is a particular case of our general setting and the most studied example in the literature. In this section, we give a full description of minimax and adaptive rates.

We check that our estimation procedure attains the minimax and adaptive rates (when known) in this context and that it provides the rates for the other setups. Very few results are available on the optimality of the rates in the adaptive setup and we prove here such results.

Let $\Lambda = [\underline{b}, \overline{b}] \times [\underline{a}, \overline{a}] \times [\underline{r}, \overline{r}] \times [\underline{L}, \overline{L}] \subset [0,\infty) \times [0,\infty) \times (0,2] \times (0,\infty)$ be a set of parameters $\lambda = (b, a, r, L)$. We shall denote by $\varphi_n$ the minimax rate of convergence over the class $\mathcal{S}(\lambda)$; see, for example, Butucea [4] for a definition. We shall say that an estimator is adaptive minimax over the family of classes $\mathcal{S}(\lambda), \lambda \in \Lambda$, if it attains the minimax rate $\varphi_n$ uniformly in $\lambda$.

It is not always possible to attain the minimax rate uniformly over a set of parameters $\Lambda$. It may happen that there is a loss in the rate due to adaptation, see Lepski [26]. We shall say that an estimator is adaptive for the adaptive rate $\phi_n$ if it attains this rate uniformly in $\lambda$ over $\Lambda$ and if, moreover, the lower bounds hold for this rate uniformly in $\lambda$ over $\Lambda$. For a definition, see Butucea [3].

For pointwise estimation of $g$, we can take $\psi(x) = \delta_{\{x_0\}}(x)$ for any given $x_0$, where $\delta_{\{x_0\}}$ is the Dirac measure at $x_0$. This implies $\psi^*(t) = \mathrm{e}^{\mathrm{i}tx_0}$ and $|\psi^*(t)| = 1$. Therefore, the rates correspond to the particular case $B = A = R = 0$ in (4) and in Table 1. They are summarized more simply in Table 2. Our procedures attain the rates already found in pointwise deconvolution and cover all other previously unknown setups.

When $r > 0, \rho > 0$, the value of $m_n$ is not explicitly given. It is obtained as the solution of the equation

$$m_n^{2b+2\gamma+(1-\rho)_+} \exp\{2\alpha(\pi m_n)^\rho + 2a(\pi m_n)^r\} = \mathrm{O}(n). \tag{20}$$

Consequently, the rate of $\hat{g}_{m_n}$ is not easy to give explicitly and depends on the ratio $r/\rho$. If $r/\rho$ or $\rho/r$ belongs to $]k/(k+1);(k+1)/(k+2)]$ with integer $k$, the rate of convergence can be expressed as a function of $k$. For explicit formulae for the rates, see Lacour [23].

These rates are known to be optimal in the minimax sense as indicated in Table 2. The case $r = 0$ is studied by Fan [16], the case $r = 0, \rho > 0$ by Cator [9] and the case $r > 0, \rho = 0$ by Butucea [4]. The rate in the case $r > 0$, $\rho > 0$, $\gamma = 0$ is proven optimal in the minimax sense in Butucea and Tsybakov [6] for $r \le \rho$. By using their construction

**Table 2.** Choice of $m_n$ for pointwise deconvolution and corresponding rates under assumptions (2) and (3). Adaptive rates for comparison. $B_m$ is abbreviated for $m^{-2b+1-r}\exp(-2a(\pi m)^r)$ and $V_m$ for $m^{2\gamma+1-\rho+(1-\rho)_+}\exp(2\alpha(\pi m)^\rho)/n$

| | $\rho=0$ ordinary smooth | $\rho>0$ supersmooth |
|---|---|---|
| $r=0$ Sob.$(b)$ | $\pi m_n = n^{1/(2b+2\gamma)}$ <br> $\varphi_n^2=\mathrm{O}(n^{-(2b-1)/(2b+2\gamma)})$ <br> *minimax rate* (Fan [16]) <br> $\phi_n^2=\mathrm{O}((n/\ln(n))^{-(2b-1)/(2b+2\gamma)})$ <br> *adaptive rate* (NEW) | $\pi m_n=[\ln(n)/(2\alpha+1)]^{1/\rho}$ <br> $\varphi_n^2=\mathrm{O}((\ln(n))^{-(2b-1)/\rho})$ <br> *minimax rate* (Fan [16]) <br> $\phi_n^2=\mathrm{O}((\ln(n))^{-(2b-1)/\rho})$ <br> *adaptive minimax rate (no loss)* (Cator [9]) |
| $r>0$ $\mathcal{C}^\infty$ | $\pi m_n=[\ln(n)/2b]^{1/r}$ <br> $\varphi_n^2=\mathrm{O}(\frac{\ln(n)^{(2\gamma+1)/r}}{n})$ <br> *minimax rate* (Butucea [3]) <br> $\phi_n^2=\mathrm{O}(\frac{\ln\ln(n)\ln(n)^{(2\gamma+1)/r}}{n})$ <br> *adaptive rate* (NEW) | $m_n$ solution of (20) <br> $=\ln(n)-(\ln\ln(n))^2$ <br> $\varphi_n^2=\mathrm{O}(B_{m_n})$ : *minimax rate if* $r<\rho$ <br> $\hookrightarrow$ (Butucea and Tsybakov [6]) <br> $\varphi_n^2=\mathrm{O}(V_{m_n})$ : *minimax rate if* $r\geq\rho$, $\rho\leq 1$(NEW) <br> $\phi_n^2=\mathrm{O}(m_n^{\rho I(r\geq\rho)}\varphi_n^2)$ <br> *adaptive minimax rate if* $r<\rho$ (no loss) <br> $\hookrightarrow$ (Butucea and Tsybakov [6]) <br> *adaptive rate if* $r\geq\rho$, $\rho\leq 1$ (NEW) |

and by following the same proof, we get near optimality (within a log factor) in the case $r>\rho$.

Very few results on adaptive pointwise estimation are available. We use $|\psi^*(x)|=1$ in the procedure described in Section 3, with $c_m\leq\sigma_m^2$ and $x_m^2c_m^2\leq Cx_m\sigma_m^2$ for all the choices of $x_m$ that will be found. Clearly, if $f_\varepsilon$ is ordinary smooth, the choice $x_m=(2\gamma+3)\ln(m)$ suits and if $f_\varepsilon$ is supersmooth, we can choose $x_m=4\alpha(\pi m)^\rho$. These choices coincide with the general case detailed above for $B=0$. Then we have $\sum_{m\in\mathcal{M}}\mathrm{e}^{-x_m}\omega_m^2\leq C/n$. This implies that

$$\mathbb{E}[(\hat\theta_{\hat m}-\theta)^2]\leq C\inf_{m\in\mathcal{M}}\left(\left(\int_{\pi m}^{+\infty}|g^*|\right)^2+\frac{x_m}{n}\min\left\{\int_{-\pi m}^{\pi m}|f_\varepsilon^*|^{-2},\left(\int_{-\pi m}^{\pi m}|f_\varepsilon^*|^{-1}\right)^2\right\}\right)+\frac{C'}{n}.$$

The rates still correspond to the particular case $B=A=R=0$ in (4) which are summarized in Table 2.

Let us mention that in the cases $\rho>0,\alpha>0$ and $r>0,a>0$ (i.e., both $f_\varepsilon$ and $g$ are supersmooth), then $x_m$ is of order $m^\rho$. There is no loss due to adaptation if $r<\rho$ as noticed earlier by Butucea and Tsybakov [6], but, surprisingly, we notice a loss of order $[\ln(n)]^{\rho/r}$ if $r>\rho$ associated to a rate faster than any power of logarithm. If $r=\rho$, the loss is logarithmic and the rate polynomial.

The previously defined estimator $\hat\theta_{m_n}$ with $m_n$ defined in Table 2 is adaptive minimax in the cases: ($r=0$ and $\rho>0$) and ($r>0$, $\rho>0$ and $r<\rho$). As we already noticed, estimators $\hat{\theta}_{\hat{m}}$, which are free of parameters, may attain a slower rate of convergence $\phi_n$, that is, it may happen that $\varphi_n = \mathrm{o}(\phi_n)$. Therefore, we check that the loss with respect to the minimax rate, when it occurs, is unavoidable.

**Theorem 4.1.** *The rates $\phi_n$ defined in Table 2 are adaptive rates and when either $\rho = 0$ or ($r \geq \rho > 0$ and $\rho \leq 1$) the loss with respect to the minimax rate which appears (compare in Table 2, $\varphi_n^2$ and $\phi_n^2$) is optimal, that is, it satisfies the following lower bounds:*

$$\inf_{\theta_n} \sup_{\lambda \in \Lambda} \sup_{g \in \mathcal{S}(\lambda)} \phi_n^{-2} \mathbb{E}_g[|\theta_n - \theta(g)|^2] \geq c$$

*for $n$ large enough, where the infimum is taken over all possible estimators $\theta_n$, under the additional hypothesis that the noise density is three-times continuously differentiable and*

$$\text{for polynomial noise } |f_\varepsilon'(u)| \leq C \frac{1}{|u|^{\gamma+1}}, \qquad \text{as } |u| \to \infty \tag{21}$$

$$\text{for exponential noise } |f_\varepsilon'(u)| \leq C|u|^{\rho-1} \exp(-\alpha |u|^\rho), \qquad \text{as } |u| \to \infty. \tag{22}$$

*Moreover, when $r > 0$, $r \geq \rho$ and $0 < \rho \leq 1$ the rate $\varphi_n^2$ is the minimax rate of estimation.*

***Remark 4.1.*** Note that the adaptive property of $\hat{\theta}_{\hat{m}}$ in the case $r \geq \rho$ is proved only for $\rho \leq 1$, which is a technical restriction. Nevertheless, it is worth noticing that, still under the restriction that $\rho \leq 1$, we obtain as a by-product in Theorem 4.1 the minimaxity of the rate for $r \geq \rho$. This is a new result since the latest result on the subject was proving minimaxity in the case $r < \rho$ only (see Butucea and Tsybakov [6]).

## 4.2. Stochastic volatility model

In this section, we consider the discrete time stochastic volatility model. Let $\eta_i$ be an i.i.d. centered noise process, $\mathbb{E}(\eta_i^2) = 1$ and let $V_i$ be a sequence of positive random variables. Assume that we observe $U_1, \dots, U_n$, where

$$U_i = \sqrt{V_i}\eta_i, \qquad i = 1, \dots, n. \tag{23}$$

Then the conditional variance of $U_i$ given $V_i$ equals $V_i$ which explains that $V_i$ is called the volatility process. In many contexts, this process is the process of interest. We assume moreover, $(V_i)$ and $(\eta_i)$ are independent and $(V_i)$ is a stationary $\beta$-mixing process with $\beta$-mixing coefficients denoted by $(\beta_k)$. When this model is obtained as the discretization of a set of continuous time stochastic differential equations, $V_i$ is indeed geometrically $\beta$-mixing, and $\eta_i \sim \mathcal{N}(0,1)$; see Comte and Genon-Catalot [13].

Model (23) is also considered in this form by van Es *et al.* [29] among others, under the assumption $\eta_i \sim \mathcal{N}(0,1)$. Setting

$$Z_i = \ln(U_i^2), \qquad X_i = \ln(V_i) \quad \text{and} \quad \varepsilon_i = \ln(\eta_i^2)$$

allows us to write (23) in the form (1). Then, we note that if $\eta_1 \sim \mathcal{N}(0,1)$,

$$f_\varepsilon^*(x) = \frac{2^{\mathrm{i}x}}{\sqrt{\pi}}\Gamma(1+\mathrm{i}x) \quad \text{and} \quad |f_\varepsilon^*(x)| \sim_{|x|\to+\infty} \sqrt{2/\mathrm{e}}\mathrm{e}^{-\pi|x|/2}, \tag{24}$$

by using the Stirling formula $\Gamma(z) \sim_{|z|\to+\infty} \sqrt{2\pi}z^{z-1/2}\mathrm{e}^{-z}$. We recognize (2) with $\gamma = 0$, $\alpha = \pi/2$ and $\rho = 1$.

Applying the results of Section 4.1 in the mixing context (D1) (see Proposition 2.2 and Corollary 3.1), we deduce that, if $V$ is geometrically $\beta$-mixing, we have a pointwise estimator of $g$,

$$\hat{g}_m(x) = \frac{1}{2\pi n}\int_{|t|\le\pi m} \frac{\mathrm{e}^{\mathrm{i}t(x+Z_k)}}{f_\varepsilon^*(t)}\,\mathrm{d}t$$

for which we can propose an automatic selection of $m$ which reaches the adaptive or adaptive minimax rate. The resulting rate is a negative power of $\ln(n)$ if $g$ is in a Sobolev space but it is much faster if $g$ is supersmooth (a case which is easy to meet; see the examples in Comte and Genon-Catalot [13]). Therefore, we recover as a particular case, and substantially improve the result of van Es *et al.* [29], who propose a non-adaptive kernel estimator of $g$, assuming that $g$ is known to be twice continuously differentiable.

Now, extensions of the class of discrete time stochastic volatility models have been studied (see Genon-Catalot and Kessler [17] or Chaleyat-Maurel and Genon-Catalot [11]) and, in particular, it is natural to consider more general types of distributions for $\eta$. For instance, we suppose now that $\eta^2$ follows a Gamma distribution, that is, $f_{\eta_1^2}(x) = (\mathrm{e}^{-x}x^{p-1}/\Gamma(p))\mathbf{I}_{x>0}$. In that case, we find

$$f_\varepsilon^*(x) = \frac{\Gamma(\mathrm{i}x+p)}{\Gamma(p)} \quad \text{and} \quad |f_\varepsilon^*(x)| \sim_{|x|\to+\infty} \frac{\sqrt{2\pi}\mathrm{e}^{-p}}{\Gamma(p)}|x|^{p-1/2}\mathrm{e}^{-\pi|x|/2}, \tag{25}$$

that is, $\varepsilon$ is supersmooth with $\gamma = p - 1/2$, $\alpha = \pi/2$ and $\rho = 1$ in (2). The Gaussian case corresponds to $p = 1/2$. Let us recall that the Laplace transform $Lu$ of a real valued function $u$ is defined by $Lu(x) = \int \mathrm{e}^{-xt}u(t)\,\mathrm{d}t$ as soon as it exists, and the Laplace transform of a non-negative random value $Y$ is defined by $\mathbb{E}(\mathrm{e}^{-\lambda Y})$. In this context, let $\pi$ denote the density of $V_1$, and consider that we are interested in estimating the Laplace transform of $V_1$. In fact, our general method provides an estimator of $h(\lambda) = -(L\pi)'(\lambda) = \mathbb{E}(V_1\mathrm{e}^{-\lambda V_1})$, that is, minus the derivative of the Laplace transform of $\pi$. In other words, we can estimate $h(\lambda) = \langle\psi_\lambda, g\rangle = \mathbb{E}(V_1\mathrm{e}^{-\lambda V_1}) = \mathbb{E}(\mathrm{e}^{X_1-\lambda\mathrm{e}^{X_1}})$. Actually we have, for $\lambda > 0$,

$$h(\lambda) = \langle\psi_\lambda, g\rangle, \qquad \text{with } \psi_\lambda(x) = \mathrm{e}^{x-\lambda\mathrm{e}^x},$$

and

$$\psi_\lambda^*(x) = \lambda^{-1-\mathrm{i}x}\Gamma(1+\mathrm{i}x) \sim_{|x|\to+\infty} \frac{\sqrt{2\pi}}{\mathrm{e}\lambda}\sqrt{|x|}\mathrm{e}^{-\pi|x|/2}, \tag{26}$$

(i.e., $B=1/2$, $A=\pi/2$ and $R=1$ in (4)). Let us define

$$\hat{h}_m(\lambda)=\frac{1}{2\pi n}\sum_{k=1}^{n}\int_{|t|\le\pi m}\mathrm{e}^{\mathrm{i}tZ_k}\frac{\psi_\lambda^*(t)}{f_\varepsilon^*(t)}\,\mathrm{d}t \tag{27}$$

with $f_\varepsilon^*$ and $\psi_\lambda^*$ given by (25) and (26). Then, taking into account the orders of $f_\varepsilon^*$ and $\psi_\lambda^*$, we obtain, by applying inequality (11) of Proposition 2.2 and if $p\neq 3/2$:

$$\mathbb{E}[(\hat{h}_m(\lambda)-h(\lambda))^2]\le Km\mathrm{e}^{-\pi^2 m}+\frac{K'm^{(3-2p)\vee 0}}{n}+\frac{K''\sum_{k\ge 0}\beta_k}{n},$$

where $K$, $K'$ and $K''$ are positive constants, $K''=2(\int|\psi^*|)^2$. If $p=3/2$, the variance term has order $\ln(m)/n$. Then notice that (D1) is satisfied in our model. Therefore, we get

**Proposition 4.1.** *Consider model (23) with* (D*1*), *(25) and (26). Assume that* $(X_k)=(\ln(V_k))$ *is* $\beta$*-mixing with* $\sum_k\beta_k<+\infty$, *then* $\hat{h}_m$ *defined by (27) satisfies, for* $\lambda>0$,

$$\mathbb{E}[(\hat{h}_m(\lambda)-h(\lambda))^2]$$
$$\le Km\mathrm{e}^{-\pi^2 m}+\frac{K'(m^{(3-2p)\vee 0}\mathbf{I}_{p\neq 3/2}+\ln(m)\mathbf{I}_{p=3/2})}{n}+\frac{K''\sum_k\beta_k}{n},$$

*where* $K$, $K'$ *and* $K''$ *are positive constants.*

In other words, using Table 1 we obtain a rate of order $[\ln(n)]^{(3-2p)\vee 1}/n$ (i.e., always less than $\ln^3(n)/n$), whatever the smoothness of $g$ is.

No adaptation is required if $p>3/2$. If $p\le 3/2$, the risk of the adaptive estimator is obtained by applying Corollary 3.1 and by choosing $x_m=4\ln(m)$:

**Proposition 4.2.** *Consider the stochastic volatility model (23) with* (D*1*), *(25) and (26). Assume that* $(X_i)$ *is geometrically* $\beta$*-mixing and consider* $\hat{h}_m$ *defined by (27), with* $\hat{m}$ *defined by (17). For any* $\lambda>0$, *and* $p\le 3/2$

$$\mathbb{E}[(\hat{h}_{\hat{m}}(\lambda)-h(\lambda))^2]$$
$$\le K\inf_{m\in\mathcal{M}}\left[\left(\int_{|u|\ge\pi m}|g^*(u)\psi_\lambda^*(u)|\,\mathrm{d}u\right)^2\right.$$
$$\left.+\frac{(m^{3-2p}\mathbf{I}_{p<3/2}+\ln(m)\mathbf{I}_{p=3/2})\ln(m)}{n}\right]+K'\frac{\ln(n)}{n}.$$

This corresponds to the case where a loss of order $\ln(\ln(n))$ occurs with respect to the non-adaptive rate.

***Remark 4.2.*** The Gaussian case, for $p = 1/2$ is not especially studied here because another strategy is available. Indeed for $\eta \sim \mathcal{N}(0,1)$, $\mathbb{E}(\mathrm{e}^{\mathrm{i}\sqrt{2\lambda}U_1}) = \mathbb{E}[\mathbb{E}(\mathrm{e}^{\mathrm{i}\sqrt{2\lambda V_1}\eta_1}|V_1)] = \mathbb{E}(\mathrm{e}^{-\lambda V_1})$. Therefore the Laplace transform of $\pi$, $L\pi(\lambda)$ can be directly estimated by an empirical mean of the $\exp(i\sqrt{2\lambda}U_k)$, which is an unbiased estimator reaching the parametric rate $1/n$. The rate would be the same for estimating $h$, as by differentiating,

$$h(\lambda) = \mathbb{E}(V_1\mathrm{e}^{-\lambda V_1}) = (-\mathrm{i}/\sqrt{2\lambda})\mathbb{E}(U_1\mathrm{e}^{\mathrm{i}\sqrt{2\lambda}U_1}).$$

The method above reaches for $p = 1/2$, the rate $\ln^w(n)\ln(\ln(n))/n$, where $1 \le w \le 2$. Therefore, it is not optimal for any $p$. But the last strategy here exploits an additional assumption ($\eta$ is Gaussian) which the general methods do not take into account.

# 5. Proofs

## 5.1. Selected $m$ for Table 1

The squared bias variance compromise is performed via the following choices of $m$, denoted by $m_n$, in the cases enumerated in Table 1:

(1) (2) and (3) (a) Case $r \vee R = 0$, $m_n = \mathrm{O}(n^{1/(2b+2B-1)})$ as $2b - 1 > 0$ when $r = 0$.
(b) Case $r \vee R > 0$, $\pi m_n = (\ln(n)/C)^{1/(r\vee R)}$ for some $C \ge A + a$.

(4) Optimal $m_n$ is such that $2a(\pi m_n)^r + 2A(\pi m_n)^R = \ln(n) - (2b + 2B - 1)\ln(m_n)$. Take,
e.g., $\pi m_n = (\ln(n)/C)^{1/(r\vee R)}$ with sufficiently large $C > 0$.

(5) Take $m_n = \mathrm{O}(n^{1/(2b+2B-1)})$.

(6) Optimal $m_n$ is such that $2a(\pi m_n)^r + 2A(\pi m_n)^R = \ln(n) - (2b + 2\gamma)\ln(m_n)$, which gives $\pi m_n = (\ln(n)/(2a) - A/a(\ln(n)/(2a))^{R/r} - (b+\gamma)/(ar)\ln\ln(n))^{1/r}$ if $r \ge R$ and
exchange $R$ and $r$ in the last expression if $R > r$. For an easier choice, take, for example,
$m_n = (\ln(n)/C)^{1/(r\vee R)}$ for $C > 0$ large enough.

(7) $m_n = \mathrm{O}(n^{1/(2b+2\gamma)})$, $b + \gamma > 0$. (8) and (9) already discussed.

(10) The optimal $m_n$ is $\pi m_n = (\ln(n)/(2\alpha) - (b+\gamma)/(\alpha\rho)\ln\ln(n))^{1/\rho}$. For a simpler form
it is sufficient to take, for example, $\pi m_n = (\ln(n)/(4\alpha))^{1/\rho}$.

The parameters $a, b, r$ of the unknown function appear several times to select $m_n$. As $g$ is unknown, and thus $a, b, r$ are unknown, it is not possible to select $m_n$ in all the cases where the rate is slower than the parametric rate $n^{-1}$.

## 5.2. Proof of Theorem 3.1

We insert here general weights $x_{j,m}$ such that

$$H(j,m) = 4\left(1+\frac{1}{\mathrm{a}}\right)(x_{j,m}\sigma_{j,m}^2 + x_{j,m}^2 c_{j,m}^2).$$

We define

$$\Xi(m) = [\theta_m - \theta(g)]^2 + \sigma_m^2 + \sup_{j\le m} x_{j,m}\sigma_{j,m}^2$$

and

$$m_{\mathrm{opt}} = \inf\left\{m\in\mathcal{M}, \mathit{Crit}(m) \le \inf_{l\in\mathcal{M}} \mathit{Crit}(l) + \frac{1}{n}\right\}.$$

It is sufficient to prove the following theorem:

**Theorem 5.1.** *There exists some positive constant $C(\mathrm{a})$ depending only on* a, *such that*

$$\mathbb{E}[(\hat\theta_{\hat m} - \theta)^2] \le C(\mathrm{a})(\mathit{Crit}(m_{\mathrm{opt}}) + \Xi(m_{\mathrm{opt}})) + C(\mathrm{a})\left(\sum_{m\in\mathcal{M}} \mathrm{e}^{-x_m}\omega_m^2 + \sum_{j\ge m_{\mathrm{opt}}} \mathrm{e}^{-x_{j,m_{\mathrm{opt}}}}\omega_{j,m}^2 + \frac{1}{n}\right),$$

*where $\omega_m^2 = \sigma_m^2 \vee c_m + 2(\sigma_m^2\vee c_m)^2$ and $\omega_{j,m}^2 = \sigma_{j,m_{\mathrm{opt}}}^2 \vee c_{j,m} + 2(\sigma_{j,m_{\mathrm{opt}}}^2 \vee c_{j,m})^2$.*

First, note that Theorem 5.1 implies Theorem 3.1. Indeed, note that for $j\ge m$, we have $\sigma_{m,j}^2 \le \sigma_j^2$ and $c_{m,j}\le c_j$. Therefore, choosing $x_{m,j} = x_j$ implies that

$$\sum_{j\ge m_{\mathrm{opt}}} \mathrm{e}^{-x_{j,m_{\mathrm{opt}}}}\omega_{j,m}^2 \le \sum_{m\in\mathcal{M}} \mathrm{e}^{-x_m}\omega_m^2.$$

Moreover $\mathit{Crit}(m) \le (\int_{|x|\ge\pi m} |\psi^*(x)g^*(x)|\,\mathrm{d}x)^2 + \mathrm{pen}(m)$ and $\Xi(m) \le (\int_{|x|\ge\pi m}|\psi^*(x)\times g^*(x)|\,\mathrm{d}x)^2 + 2\mathrm{pen}(m)$. This implies Theorem 3.1.

Now we establish the following lemma:

**Lemma 5.1.** *For all $m\in\mathcal{M} := \{1,\dots,m_n\}$, for all $x>0$,*

$$\mathbb{P}\left(\widehat{\mathit{Crit}}(m) > (1+\mathrm{a})\mathit{Crit}(m) + 4\left(1+\frac{1}{\mathrm{a}}\right)\left(x+x^2\right)\right) \le \sum_{j\ge m, j\in\mathcal{M}} \mathrm{e}^{-x_{j,m}}\mathrm{e}^{-x/(\sigma_{j,m}^2\vee c_{j,m})}.$$

**Proof.** Recall that the Bernstein inequality for a sum $S_n = \sum_{k=1}^n Y_k$ of i.i.d. random variables $Y_k$ having $\mathrm{var}(Y_1) \le v^2$ and $\|Y_1\|_\infty \le 1/\mathrm{a}$ states that

$$\mathbb{P}\left((S_n - \mathbb{E}(S_n))/n \ge \sqrt{2uv^2/n} + \frac{u}{\mathrm{a}n}\right) \le \exp(-u).$$

We put for $j \ge m$

$$Y_k = Y_k(j,m) = \frac{1}{2\pi}\int_{\pi m \le |t| \le \pi j} \mathrm{e}^{\mathrm{i}tZ_k}\frac{\psi^*(t)}{f_\varepsilon^*(t)}\,\mathrm{d}t. \tag{28}$$

Then $S_n/n = \hat\theta_j - \hat\theta_m$ and $\mathbb{E}(S_n/n) = \mathbb{E}(\hat\theta_j - \hat\theta_m) = \theta_j - \theta_m$. Moreover, we obtain that $v^2/n \le \sigma_{j,m}^2$ and $1/(\mathrm{a}n) = c_{j,m}$. It follows that

$$\mathbb{P}\{[(\hat\theta_j - \hat\theta_m) - (\theta_j - \theta_m)]^2 \ge (\sigma_{j,m}\sqrt{2u} + c_{j,m}u)^2\} \le 2\mathrm{e}^{-u}.$$

Now, from the simple fact that $(x+y)^2 \le (1+1/\mathrm{a})x^2 + (1+\mathrm{a})y^2$ for any real numbers $x$, $y$, we deduce by setting $u = y$ and $v = x + y$ that $(v-u)^2 \ge (1/(1+1/\mathrm{a}))v^2 - (1+\mathrm{a})/(1+1/\mathrm{a})u^2$. Use also the fact that $(A+B)^2 \le 2(A^2+B^2)$ for any real numbers $A$, $B$, to obtain

$$\mathbb{P}\{(\hat\theta_j - \hat\theta_m)^2 \ge (1+\mathrm{a})(\theta_j - \theta_m)^2 + 2(1+1/\mathrm{a})(2\sigma_{j,m}^2 u + c_{j,m}^2 u^2)\} \le 2\mathrm{e}^{-u}.$$

Now we set $u = x_{j,m} + x/(\sigma_{j,m}^2 \vee c_{j,m})$ and we find

$$\begin{aligned}&\mathbb{P}\left\{(\hat\theta_j - \hat\theta_m)^2 - H(j,m) \ge (1+\mathrm{a})(\theta_j - \theta_m)^2 + 4\left(1 + \frac{1}{\mathrm{a}}\right)(x + x^2)\right\}\\ &\quad\le 2\mathrm{e}^{-x_{j,m}}\mathrm{e}^{-x/(\sigma_{j,m}^2 \vee c_{j,m})}.\end{aligned}$$

To conclude we write

$$\begin{aligned}&\mathbb{P}\left(\widehat{Crit}(m) > (1+\mathrm{a})Crit(m) + 4\left(1+\frac{1}{\mathrm{a}}\right)(x+x^2)\right)\\ &\quad\le \mathbb{P}\left\{\exists j \ge m, j \in \mathcal{M}, (\hat\theta_j - \hat\theta_m)^2 - H(j,m) \ge (1+\mathrm{a})(\theta_j - \theta_m)^2 + 4\left(1+\frac{1}{\mathrm{a}}\right)(x+x^2)\right\}\\ &\quad\le 2\sum_{j\ge m, j\in\mathcal{M}} \mathrm{e}^{-x_{j,m}}\mathrm{e}^{-x/(\sigma_{j,m}^2 \vee c_{j,m})}.\end{aligned}$$

This ends the proof of Lemma 5.1. □

Now we follow the steps of the proof of Laurent *et al.* [24].

- We first consider the case where $\hat m \le m_{\mathrm{opt}}$. Following the same lines of proof, we get

$$\mathbb{P}\left(\frac{1}{2}(\hat\theta_{\hat m} - \theta(g))^2 > (1+\mathrm{a})\,Crit(m_{\mathrm{opt}}) + 4\left(1 + \frac{1}{\mathrm{a}}\right)(x + x^2)\right.$$

$$+\sup_{j\le m_{\text{opt}}} H(m_{\text{opt}},j)+(\hat\theta_{m_{\text{opt}}}-\theta(g))^2+\frac{1}{n}\cap\{\hat m\le m_{\text{opt}}\}\Big) \qquad (29)$$

$$\le \sum_{j\ge m_{\text{opt}}} \mathrm{e}^{-x_{j,m_{\text{opt}}}}\mathrm{e}^{-x/(\sigma^2 j,m_{\text{opt}}\vee c_{j,m_{\text{opt}}})}.$$

- Now we consider the case $\hat m > m_{\text{opt}}$. We apply the Bernstein inequality to

$$\tilde Y_k=\tilde Y_k(m)=\frac{1}{2\pi}\int_{|t|\le \pi m}\mathrm{e}^{\mathrm{i}tZ_k}\frac{\psi^*(t)}{f_\varepsilon^*(t)}\,\mathrm{d}t,$$

in the same way as in Lemma 5.1. We obtain, for all $m\in\mathcal{M}$,

$$\mathbb{P}\bigg((\hat\theta_m-\theta(g))^2\ge(1+\mathrm{a})(\theta_m-\theta(g))^2+4\Big(1+\frac{1}{\mathrm{a}}\Big)(x+x^2)+\mathrm{pen}(m)\bigg)$$
$$\le 2\mathrm{e}^{-x_m}\mathrm{e}^{-x/(\sigma_m^2\vee c_m)}.$$

This implies that

$$\mathbb{P}\bigg((\hat\theta_{\hat m}-\theta(g))^2\ge(1+\mathrm{a})(\theta(g_{\hat m})-\theta(g))^2+4\Big(1+\frac{1}{\mathrm{a}}\Big)(x+x^2)+\mathrm{pen}(\hat m)\bigg)$$
$$\le \sum_{m\in\mathcal{M}}2\mathrm{e}^{-x_m}\mathrm{e}^{-x/(\sigma_m^2\vee c_m)}.$$

As $\sup_{j\ge m}[(\hat\theta_m-\hat\theta_j)^2-H(j,m)]\ge(\hat\theta_m-\hat\theta_m)^2-H(m,m)=0$, we have $\widehat{Crit}(m)\ge \mathrm{pen}(m)$. Using the inequalities $\mathrm{pen}(m)\le\widehat{Crit}(\hat m)\le\widehat{Crit}(m_{\text{opt}})+1/n$, we obtain

$$\mathbb{P}\bigg((\hat\theta_{\hat m}-\theta(g))^2\ge(1+\mathrm{a})(\theta_{\hat m}-\theta(g))^2+4\Big(1+\frac{1}{\mathrm{a}}\Big)(x+x^2)+\widehat{Crit}(m_{\text{opt}})+\frac{1}{n}\bigg)$$
$$\le \sum_{m\in\mathcal{M}}2\mathrm{e}^{-x_m}\mathrm{e}^{-x/(\sigma_m^2\vee c_m)}.$$

If $\hat m>m_{\text{opt}}$, then $(\hat\theta_m-\theta(g))^2\le\sup_{j\ge m_{\text{opt}}}(\theta_j-\theta(g))^2$ and we apply Lemma 5.1 with $m=m_{\text{opt}}$. This yields

$$\mathbb{P}\bigg((\hat\theta_{\hat m}-\theta(g))^2\ge(1+\mathrm{a})\bigg(\sup_{j\ge m_{\text{opt}}}(\theta_j-\theta(g))^2+8\Big(1+\frac{1}{\mathrm{a}}\Big)(x+x^2)$$
$$+(1+\mathrm{a})\,Crit(m_{\text{opt}})+\frac{1}{n}\cap\{\hat m>m_{\text{opt}}\}\bigg)\bigg)$$
$$\le\sum_{m\in\mathcal{M}}2\mathrm{e}^{-x_m}\mathrm{e}^{-x/(\sigma_m^2\vee c_m)}+\sum_{j\ge m_{\text{opt}}}2\mathrm{e}^{-x_{j,m_{\text{opt}}}}\mathrm{e}^{-x/(\sigma_{j,m_{\text{opt}}}^2\vee c_{j,m_{\text{opt}}})}. \qquad (30)$$

Let

$$C_{m_{\rm opt}} = 3(1+{\rm a})\,Crit(m_{\rm opt}) + 2\sup_{j\le m_{\rm opt}} H(m_{\rm opt}, j) + (1+{\rm a})\sup_{j\ge m_{\rm opt}} (\theta_j - \theta(g))^2 + \frac{3}{n}$$

and

$$X = (\hat\theta_{\hat m} - \theta(g))^2, \qquad Y = 2(\hat\theta_{m_{\rm opt}} - \theta(g))^2.$$

It follows from (29) and (30) that, for all $x>0$,

$$\begin{aligned}&\mathbb{P}\biggl(X - Y > C_{m_{\rm opt}} + 24\Bigl(1+\frac{1}{\rm a}\Bigr)(x\vee x^2)\biggr)\\ &\quad\le \sum_{m\in\mathcal{M}} 2{\rm e}^{-x_m}{\rm e}^{-x/(\sigma_m^2\vee c_m)} + \sum_{j\ge m_{\rm opt}} 2{\rm e}^{-x_{j,m_{\rm opt}}}{\rm e}^{-x/(\sigma_{j,m_{\rm opt}}^2\vee c_{j,m_{\rm opt}})}.\end{aligned}$$

We write that $\mathbb{E}(X) = \mathbb{E}(X\mathbf{I}_{X\ge Y+C_{m_{\rm opt}}}) + \mathbb{E}(X\mathbf{I}_{X\le Y+C_{m_{\rm opt}}}) \le \mathbb{E}[(X-Y-C_{m_{\rm opt}})_+] + \mathbb{E}(Y+C_{m_{\rm opt}})$.

Then, setting $C_{\rm a} = 24(1+1/{\rm a})$ and $Z = X - Y - C_{m_{\rm opt}}$

$$\begin{aligned}\mathbb{E}[Z_+] &= \int_0^{+\infty}\mathbb{P}(Z>t)\,{\rm d}t = C_{\rm a}\biggl(\int_0^1 \mathbb{P}(Z>C_{\rm a}u)\,{\rm d}u + \int_1^\infty \mathbb{P}(Z>C_{\rm a}u)\,{\rm d}u\biggr)\\ &= C_{\rm a}\biggl(\int_0^1 \mathbb{P}(Z>C_{\rm a}(u\vee u^2))\,{\rm d}u + 2\int_1^\infty \mathbb{P}(Z>C_{\rm a}(v\vee v^2))v\,{\rm d}v\biggr),\end{aligned}$$

$$\begin{aligned}\mathbb{E}[(X-Y-C_{m_{\rm opt}})_+] &\le C_{\rm a}\sum_{m\in\mathcal{M}} 2{\rm e}^{-x_m}(\sigma_m^2\vee c_m + 2(\sigma_m^2\vee c_m)^2)\\ &\quad + C_{\rm a}\sum_{j\ge m_{\rm opt}} 2{\rm e}^{-x_{j,m_{\rm opt}}}(\sigma_{j,m_{\rm opt}}^2\vee c_{j,m_{\rm opt}} + 2(\sigma_{j,m_{\rm opt}}^2\vee c_{j,m_{\rm opt}})^2)\\ &= C_{\rm a}\biggl(\sum_{m\in\mathcal{M}} 2{\rm e}^{-x_m}\omega_m^2 + \sum_{j\ge m_{\rm opt}} 2{\rm e}^{-x_{j,m_{\rm opt}}}\omega_{j,m_{\rm opt}}^2\biggr).\end{aligned}$$

The end of the proof is the same as in Laurent *et al.* [24].

## 5.3. Proof of Proposition 2.2

The same decomposition of the risk and upper bound for the bias hold, as in Section 2.2. Only the variance has to be re-examined. The basic idea is that, for $k\ne\ell$, ${\rm cov}({\rm e}^{{\rm i}tZ_k}, {\rm e}^{{\rm i}sZ_\ell}) = f_\varepsilon^*(t)f_\varepsilon^*(-s){\rm cov}({\rm e}^{{\rm i}tX_k}, {\rm e}^{{\rm i}sX_\ell})$ by conditioning on $(X_k, X_\ell)$. The additional trick is the standard covariance inequality for $\beta$-mixing variables (see, e.g.,

Doukhan [15]), which implies that $|\mathrm{cov}(\mathrm{e}^{\mathrm{i}tX_k}, \mathrm{e}^{\mathrm{i}sX_\ell})| \le \beta_{|k-\ell|}$.

$$\begin{aligned}\mathrm{Var}(\hat{\theta}_m) &= \frac{1}{4\pi^2 n^2} \sum_{k,\ell=1,k\neq\ell}^{n} \int_{-\pi m}^{\pi m}\int_{-\pi m}^{\pi m} \mathrm{cov}(\mathrm{e}^{\mathrm{i}tX_k}, \mathrm{e}^{\mathrm{i}sX_\ell})\psi^*(t)\psi^*(-s)\,\mathrm{d}s\,\mathrm{d}t \\ &\quad + \frac{1}{4\pi^2 n^2} \sum_{k=1}^{n} \int_{-\pi m}^{\pi m}\int_{-\pi m}^{\pi m} \mathrm{cov}(\mathrm{e}^{\mathrm{i}tZ_k}, \mathrm{e}^{\mathrm{i}sZ_k})\frac{\psi^*(t)}{f_\varepsilon^*(t)}\frac{\psi^*(-s)}{f_\varepsilon^*(-s)}\,\mathrm{d}s\,\mathrm{d}t.\end{aligned} \tag{31}$$

The last term is the standard variance term of the independent case. The first one is bounded in modulus by

$$\begin{aligned}&\frac{2}{4\pi^2 n^2} \sum_{1\le k<\ell\le n}^{n} \int_{-\pi m}^{\pi m}\int_{-\pi m}^{\pi m} |\mathrm{cov}(\mathrm{e}^{\mathrm{i}tX_1}, \mathrm{e}^{\mathrm{i}sX_{\ell-k}})||\psi^*(t)\psi^*(-s)|\,\mathrm{d}s\,\mathrm{d}t \\ &\le \frac{1}{2\pi n}\sum_{k=1}^{n}\beta_k\left(\int_{-\pi m}^{\pi m}|\psi^*(t)|\,\mathrm{d}t\right)^2.\end{aligned}$$

This gives the result.

## 5.4. Proof of Proposition 2.3

Under (D2), we only obtain that for $k<\ell$, $\mathrm{cov}(\mathrm{e}^{\mathrm{i}tZ_k}, \mathrm{e}^{\mathrm{i}sZ_\ell}) = f_\varepsilon^*(-s)\,\mathrm{cov}(\mathrm{e}^{\mathrm{i}tZ_k}, \mathrm{e}^{\mathrm{i}sX_\ell})$ by conditioning on $(X_\ell)$. The covariance inequality for $\beta$-mixing variables (see, e.g., Doukhan [15]) still applies (but to the variables $(X_k, Z_k)$ and $(X_\ell, Z_\ell)$ and implies that $|\mathrm{cov}(\mathrm{e}^{\mathrm{i}tZ_k}, \mathrm{e}^{\mathrm{i}sX_\ell})| \le \beta_{|k-\ell|}$. Then (31) remains true but leads, for the bound of the modulus of the last term, to:

$$\begin{aligned}&\frac{2}{4\pi^2 n^2} \sum_{1\le k<\ell\le n}^{n} \int_{-\pi m}^{\pi m}\int_{-\pi m}^{\pi m} |\mathrm{cov}(\mathrm{e}^{\mathrm{i}tZ_1}, \mathrm{e}^{\mathrm{i}sX_{\ell-k}})|\left|\frac{\psi^*(t)}{f_\varepsilon^*(t)}\psi^*(-s)\right|\mathrm{d}s\,\mathrm{d}t \\ &\le \frac{1}{2\pi n}\sum_{k=1}^{n}\beta_k\left(\int_{-\pi m}^{\pi m}|\psi^*(t)|\,\mathrm{d}t\right)\left(\int_{-\pi m}^{\pi m}\left|\frac{\psi^*(t)}{f_\varepsilon^*(t)}\right|\mathrm{d}t\right).\end{aligned}$$

This gives inequality (14).

For the proof of (16), the result follows from the inequality

$$\left(\int_{-\pi m}^{\pi m}|\psi^*|(t)\,\mathrm{d}t\right)\left(\int_{-\pi m}^{\pi m}|\psi^*/f_\varepsilon^*|(t)\,\mathrm{d}t\right) \le \left(\int_{-\pi m}^{\pi m}|\psi^*/f_\varepsilon^*|(t)\,\mathrm{d}t\right)^2$$

and the fact that the new mixing term is always negligible with respect to the independent variance term if $\varepsilon$ is supersmooth (case $A, \rho > 0$). If $\varepsilon$ is ordinary smooth, then we only have to study when $m^{(-B+1)_+ + \gamma - B + 1}$ is less than $m^{2\gamma - 2B + 1}$, which occurs as soon as $\gamma > \max(B, 1)$.

### 5.5. Proof of Corollary 3.1

The main difference with respect to the proof of Theorem 3.1 lies in the Bernstein inequality which must be written in the mixing context. For geometrically mixing variables (and $q = q_n = 2\ln(n)/c$ if $\beta_k \leq \mathrm{e}^{-ck}$), we get from Theorem 4, page 36 in Doukhan [15] that

$$\mathbb{P}\left(\frac{S_n - \mathbb{E}(S_n)}{n} \geq \sqrt{\frac{2u\tilde{v}^2}{n}} + \frac{2\ln(n)u}{can}\right) \leq \mathrm{e}^{-u} + \frac{2}{n^2},$$

with $\|Y_1\|_\infty \leq 1/\mathrm{a}$ and $(1/q)\mathrm{Var}\left(\sum_{k=1}^q Y_k\right) \leq \tilde{v}^2$.

In all cases, $|\mathcal{M}| \leq n$, so that summing up the residuals of order $1/n^2$ will give negligible terms of order $1/n$. Next, the variables are still given by (28) with $c_{j,m}$ and $c_m$ the same as previously multiplied by $2\ln(n)/c$. This gives $\tilde{c}_{j,m} = (2\ln(n)/2)c_{j,m}$ and $\tilde{c}_m = (2\ln(n)/c)c_m$. At last, it follows from the above computation of $\mathrm{Var}(\breve{\theta}_m)$ that the new variance terms denoted by $\tilde{\sigma}^2_{j,m}$, $\tilde{\sigma}^2_m$ can be bounded under (D1) by

$$\tilde{\sigma}^2_{j,m} \leq \sigma^2_{j,m} + \frac{1}{\pi n}\sum_{k\geq 1}\beta_k\left(\int_{\pi(m\wedge j)\leq |t| \leq \pi(m\vee j)} |\psi^*(t)|\,\mathrm{d}t\right)^2,$$

and analogously for $\tilde{\sigma}^2_m$. It follows from our set of assumptions that $\tilde{\sigma}^2_{j,m} \leq \sigma^2_{j,m} + c/n \leq 2\sigma^2_{j,m}$ and $\tilde{\sigma}^2_m \leq 2\sigma^2_m$. The case (D2) is analogous under the given more restrictive assumptions. The Corollary 3.1 follows.

### 5.6. Proof of Theorem 4.1

We describe first the general procedure for proving the theorem and postpone details of constructions and proofs to Section 5.6. As the adaptation loss is different according to whether $r = 0$ or $r \neq 0$, respectively $\rho = 0$ or $\rho \neq 0$, explicit constructions are needed for each of the following setups: (1) $r = 0$, $\rho = 0$; (2) $b = 0$, $r > 0$, $\rho = 0$; (3) $b = 0$, $r > 0$, $0 < \rho \leq 1$ and $r \geq \rho$. We take $b = 0$ without loss of generality, in order to simplify polynomial factors in our explicit constructions.

Typically, we construct two probability densities $g_0 \in \mathcal{S}(\overline{\lambda})$ and $g_{1,n} \in \mathcal{S}(\underline{\lambda})$ where $\overline{\lambda}, \underline{\lambda} \in \Lambda$. Moreover

$$g_{1,n}(x) = g_0(x) + G(x - x_0, m) \qquad \text{for } m = m_n \to \infty \text{ with } n \text{ and } \int G(\cdot, m) = 0 \ \forall m.$$

Note that the likelihoods of the model become $f_0^Z = g_0 \star f_\varepsilon$ under $g_0$ and

$$f_{1,n}^Z(x) = [g_{1,n} \star f_\varepsilon](x) = f_0^Z(x) + [G(\cdot, m) \star f_\varepsilon](x - x_0)$$

under $g_{1,n}$. Then

$$\inf_{\theta_n} \sup_{\lambda\in\Lambda} \sup_{g\in\mathcal{S}(\lambda)} \phi_{n,\lambda}^{-2}\mathbb{E}_g[|\theta_n - \theta(g)|^2]$$

$$\geq \inf_{\theta_n} \max\{\phi_{n,\overline{\lambda}}^{-2}\mathbb{E}_{g_0}[|\theta_n-\theta(g_0)|^2], \phi_{n,\underline{\lambda}}^{-2}\mathbb{E}_{g_{1,n}}[|\theta_n-\theta(g_{1,n})|^2]\}$$
$$\geq \inf_{T_n} \max\{q_n^2\mathbb{E}_{g_0}[T_n^2], \mathbb{E}_{g_{1,n}}[|T_n-G(0,m)/\phi_{n,\underline{\lambda}}|^2]\},$$

where $q_n=\phi_{n,\underline{\lambda}}/\phi_{n,\overline{\lambda}}\to\infty$ when $n\to\infty$, with a proper choice of $\underline{\lambda}$, $\overline{\lambda}$ and $T_n=(\theta_n-\theta(g_0))/\phi_{n,\underline{\lambda}}$.

From now on we denote $\mathbb{P}_0=\mathbb{P}_{g_0}$, $\mathbb{E}_0=\mathbb{E}_{g_0}$ and $\mathbb{P}_1=\mathbb{P}_{g_{1,n}}$, $\mathbb{E}_1=\mathbb{E}_{g_{1,n}}$. Following Theorem 6 in Tsybakov [28] we can deduce that, if $|G(0,m)/\phi_{n,\underline{\lambda}}|\geq c>0$ and if for some fixed $0<\epsilon<1$ and $\tau>0$

$$\mathbb{P}_1\left(\frac{\mathrm{d}\mathbb{P}_0}{\mathrm{d}\mathbb{P}_1}\geq\tau\right)\geq 1-\epsilon, \tag{32}$$

then

$$\inf_{T_n}\max\{q_n^2\mathbb{E}_0[T_n^2],\mathbb{E}_1[|T_n-G(0,m)/\phi_{n,\underline{\lambda}}|^2]\}\geq\frac{\tau q_n^2\epsilon^2c^4(1-\epsilon)^2}{\tau q_n^2\epsilon^2c^2+(1-\epsilon)^2c^2}. \tag{33}$$

If we can choose $\tau=\tau_n$ such that $\tau_n q_n^2\to\infty$ with $n$, then the bound from below in (33) tends to $c^2(1-\epsilon)^2$ so it will be larger than $c^2(1-\epsilon)^4>0$ for $n$ large enough. Note also that (33) may provide the exact asymptotic constant in case $c\to1$ and $\mathbb{P}_1(\mathrm{d}\mathbb{P}_0/\mathrm{d}\mathbb{P}_1\geq\tau_n)\to1$ as $n\to\infty$.

In order to deal with (32), we proceed as follows:

$$\mathbb{P}_1\left(\frac{\mathrm{d}\mathbb{P}_0}{\mathrm{d}\mathbb{P}_1}\geq\tau\right)=\mathbb{P}_1\left(\prod_{i=1}^n\frac{g_0\star f_\varepsilon}{g_{1,n}\star f_\varepsilon}(Y_i)\geq\tau\right)=\mathbb{P}_1\left(\sum_{i=1}^n\ln\left(1-\frac{G(\cdot-x_0)\star f_\varepsilon}{g_{1,n}\star f_\varepsilon}(Y_i)\right)\geq\ln(\tau)\right)$$
$$=\mathbb{P}_1\left(\frac{\sum_{i=1}^n Z_{i,n}-n\mathbb{E}_1(Z_{1,n})}{(n\mathrm{Var}_1(Z_{1,n}))^{1/2}}\geq\frac{\ln(\tau)-n\mathbb{E}_1(Z_{1,n})}{(n\mathrm{Var}_1(Z_{1,n}))^{1/2}}\right),$$

where $Z_{i,n}=\ln(1-[G(\cdot-x_0)\star f_\varepsilon](Y_i)/g_{1,n}\star f_\varepsilon(Y_i))$ form a triangular array of independent variables. Denote

$$U_n:=\frac{\sum_{i=1}^n Z_{i,n}-nE_1(Z_{1,n})}{(n\mathrm{Var}_1(Z_{1,n}))^{1/2}}.$$

We shall prove, for each setup, Lyapunov's central limit theorem for $U_n$. Moreover, we give a lower bound $\mathbb{E}_1(Z_{1,n})\geq -c_e\kappa_n$ and an upper bound for $\mathrm{Var}_1(Z_{1,n})\leq c_v\kappa_n$, where $\kappa_n$ is such that

$$\chi^2(g_0\star f_\varepsilon,g_{1,n}\star f_\varepsilon):=\int\frac{(g_{1,n}\star f_\varepsilon-g_0\star f_\varepsilon)^2}{g_{1,n}\star f_\varepsilon}\leq\kappa_n$$

as $n\to\infty$. Choose then $\tau_n\to0$ such that

$$u_n:=\frac{\ln(\tau_n)+c_en\kappa_n}{(c_vn\kappa_n)^{1/2}}\to-\infty$$

with $n$, giving that $\mathbb{P}_1(U_n \geq u_n) \geq 1-\epsilon$ for some $0<\epsilon<1$ and large enough $n$ and thus concluding the proof of the theorem.

Now, we study in more detail the different cases.

(1) Case $r=0$, $\rho=0$ and $\Lambda=[\underline{b},\overline{b}]\times[\underline{L},\overline{L}]\subset(1/2,\infty)\times(0,\infty)$.

Let us choose $g_0$ in the class $\mathcal{S}(\overline{b},\overline{L}/2)$ such that $g_0>0$ and $g_0(x)\geq c|x|^{-2}$ as $|x|\to\infty$. We choose next the function $G$ such that $G(x,m)=m^{-\underline{b}+1/2}G(mx)$ and with $G^*$ at least three-times continuously differentiable having the property

$$\frac{I(1/2\leq|u|\leq 3/4)}{c(1+u^{2\underline{b}})}\leq G^*(u)\leq\frac{I(1/4\leq|u|\leq 1)}{c(1+u^{2\underline{b}})}.$$

Here, $m=(c_0\ln(n)/n)^{-1/(2\underline{b}+2\gamma)}$. Note that $G^*(0)=\int G=0$. First, $g_{1,n}$ is a positive function with an integral equal to 1 and it belongs to $\mathcal{S}(\underline{b},\overline{L})$. Indeed, for each fixed $x$ we have $G(x,m)\to 0$ when $n\to\infty$ and as $G^*$ is three times continuously differentiable that means $|G(x,m)|\leq \mathrm{O}(|x|^{-3})=\mathrm{o}(g_0(x))$ as $|x|\to\infty$, giving that $g_{1,n}\geq 0$ for $n$ large enough. Moreover,

$$\begin{aligned}\left(\int|g_{1,n}^*(u)|^2|u|^{2\underline{b}}\,\mathrm{d}u\right)^{1/2} &\leq\left(\int|g_0^*(u)|^2|u|^{2\underline{b}}\,\mathrm{d}u\right)^{1/2}\\ &\quad+m^{-\underline{b}-1/2}\left(\int_{1/4\leq|u|/m\leq 1}|G^*(u/m)|^2|u|^{2\underline{b}}\,\mathrm{d}u\right)^{1/2}\\ &\leq\sqrt{2\pi\overline{L}/2}+\frac{C}{c}\left(\int_{1/4}^{1}\frac{|u|^{2\underline{b}}}{(1+u^{2\underline{b}})^2}\,\mathrm{d}u\right)^{1/2}\leq(2\pi\overline{L})^{1/2},\end{aligned}$$

for $c>0$ large enough. Second,

$$\left|\frac{G(0,m)}{\phi_{n,\underline{b}}}\right|=(\phi_{n,\underline{b}})^{-1}m^{-\underline{b}+1/2}\frac{1}{2\pi}\int G^*(u)\,\mathrm{d}u\geq\frac{c_0^{-\underline{b}+1/2}}{2\pi}\int_{1/2}^{3/4}\mathrm{d}u\geq c_1\cdot c_0^{-\underline{b}+1/2}>0.$$

We shall prove that (32) holds with $\tau=n^{-(2\gamma+1)/(2\underline{b}+2\gamma)}$ and together with the fact that

$$\begin{aligned}\tau q_n^2&=\tau\frac{\phi_{n,\underline{b}}^2}{\phi_{n,\overline{b}}^2}=\tau\left(\frac{\ln(n)}{n}\right)^{-(2\gamma+1)(\overline{b}-\underline{b})/((2\underline{b}+2\gamma)(2\overline{b}+2\gamma))}\\ &=(\ln(n))^{-(2\gamma+1)(\overline{b}-\underline{b})/((2\underline{b}+2\gamma)(2\overline{b}+2\gamma))}n^{(2\gamma+1)/(2\underline{b}+2\gamma)(\underline{b}+\gamma)/(\overline{b}+\gamma)}\end{aligned}$$

tends to infinity, with $n$, the proof of (33) and hence of the theorem is finished.

We can prove that for each $x_0$

$$\sup_x\frac{|[G(m(\cdot-x_0))\star f_\varepsilon](x)|}{f_0^Z(x)}=\mathrm{o}(1),\qquad\text{as } n\to\infty,\tag{34}$$

therefore $f_{1,n}^Z(x) = f_0^Z(x)(1 + \mathrm{o}(1))$, where $\mathrm{o}(1) \to 0$, $n \to \infty$ uniformly in $x$. As we chose $g > 0$ then $f_0^Z > 0$ and together with the previous statement it means that for any $M > 0$ we can find a constant $c_2 > 0$ such that $f_{1,n}^Z \geq 1/c_2$ on $[-M, M]$. Moreover, for some $M > 0$ large enough, see Butucea and Tsybakov [6], $f_0^Z(x) = g_0 \star f_\varepsilon(x) \geq C_2/x^2$, as $|x| \geq M$.

Therefore, for large enough $M > 0$, $f_{1,n}^Z(x) \geq 1/(c_3|x|^2)$, for some constant $c_3 > 0$ and for $|x| \geq M$. Finally, we deal with

$$\begin{aligned}\chi^2(f_0^Z, f_{1,n}^Z) &= m^{-2\underline{b}+1} \int \frac{[G(m(\cdot - x_0)) \star f_\varepsilon]^2(x)}{f_{1,n}^Z(x)} \mathrm{d}x \\ &\leq m^{-2\underline{b}+1} \left( c_2 \int_{|x| \leq M} [G(m(\cdot - x_0)) \star f_\varepsilon]^2(x) \mathrm{d}x \right. \\ &\qquad \left. + c_3 \int_{|x| > M} |x|^2 [G(m(\cdot - x_0)) \star f_\varepsilon]^2(x)\, \mathrm{d}x \right),\end{aligned}$$

say $T_1$ and $T_2$, for some fixed, large $M > 0$. Then

$$\begin{aligned}T_1 &\leq m^{-2\underline{b}-1} \frac{c_2}{2\pi} \int \left| G^*\left(\frac{u}{m}\right) f_\varepsilon^*(u) \right|^2 \mathrm{d}u \\ &\leq c_4 m^{-2\underline{b}-1} \int_{m/4}^{m} \frac{1}{|u|^{2\gamma}} \mathrm{d}u \leq c_5 m^{-2\underline{b}-2\gamma} \leq c_6 \frac{c_0 \ln(n)}{n}.\end{aligned} \tag{35}$$

For $T_2$ we follow the similar proof in Butucea and Tsybakov [6] and use condition (21) to get

$$\begin{aligned}T_2 &\leq m^{-2\underline{b}+1} \frac{c_3}{2\pi} \int \left| \frac{\partial}{\partial u} \left( \frac{1}{m} G^*\left(\frac{u}{m}\right) f_\varepsilon^*(u) \right) \right|^2 \mathrm{d}u \\ &\leq c_6 m^{-2\underline{b}-1} m^{-2\gamma} = \mathrm{o}(T_1), \qquad n \to \infty.\end{aligned} \tag{36}$$

Therefore, from (35) and (36) we have $\chi^2(f_0^Z, f_{1,n}^Z) \leq \kappa_n$, with $\kappa_n = c_\chi c_0 \ln(n)/n$. We use the fact that $-u(1+u) \leq \ln(1-u) \leq -u$ for all $u \in [0, 1/2]$ and that (34) implies that $|u| = |[G(m(\cdot - x_0)) \star f_\varepsilon](x)|/f_{1,n}^Z(x) \leq 1/2$ for $n$ large enough to get

$$\begin{aligned}\mathbb{E}_1[Z_{1,n}] &= \int \ln\left(1 - \frac{[G(\cdot, m) \star f_\varepsilon](x - x_0)}{f_{1,n}^Z(x)}\right) f_{1,n}^Z(x)\, \mathrm{d}x \\ &\geq -\int [G(\cdot, m) \star f_\varepsilon](x - x_0)\, \mathrm{d}x - \int \frac{[G(\cdot, m) \star f_\varepsilon]^2(x - x_0)}{f_{1,n}^Z(x)} \mathrm{d}x \\ &\geq -\chi^2(f_0^Z, f_{1,n}^Z) \geq -\kappa_n,\end{aligned}$$

for $n$ large enough. Indeed, note that $\int G(\cdot,m)=0$ and therefore $\int [G(\cdot,m)\star f_\varepsilon](x-x_0)\,\mathrm{d}x=0$. Moreover,

$$\begin{aligned}\mathrm{Var}_1(Z_{1,n}) &\le \mathbb{E}_1(Z_{1,n}^2)=\int \ln^2\left(1-\frac{[G(\cdot,m)\star f_\varepsilon](x-x_0)}{f^Z_{1,n}(x)}\right)f^Z_{1,n}(x)\,\mathrm{d}x\\ &\le \int \frac{[G(\cdot,m)\star f_\varepsilon]^2(x-x_0)}{f^Z_{1,n}(x)^2}\left(1+\frac{[G(\cdot,m)\star f_\varepsilon]^2(x-x_0)}{f^Z_{1,n}(x)}\right)^2 f^Z_{1,n}(x)\,\mathrm{d}x\\ &\le c_v\chi^2(f^Z_0,f^Z_{1,n})\le c_v\kappa_n,\end{aligned}$$

as by (34): $\sup_x|f^Z_0(x)/f^Z_{1,n}(x)|$ is bounded from above by some constant depending only on $g_0$ and $f_\varepsilon$. By similar calculations, we also check that

$$\begin{aligned}\mathrm{Var}_1(Z_{1,n}) &\ge \frac{1}{2}\mathbb{E}_1(Z_{1,n}^2)=\frac{1}{2}\int \ln^2\left(1-\frac{[G(\cdot,m)\star f_\varepsilon](x-x_0)}{f^Z_{1,n}(x)}\right)f^Z_{1,n}(x)\,\mathrm{d}x\\ &\ge \frac{1}{2}\int \frac{[G(\cdot,m)\star f_\varepsilon]^2(x-x_0)}{f^Z_{1,n}(x)}\,\mathrm{d}x\\ &\ge \frac{1}{2\|f^Z_{1,n}\|_\infty}\int [G(\cdot,m)\star f_\varepsilon]^2(x-x_0)\,\mathrm{d}x\ge c'_v\kappa_n\end{aligned}$$

and that

$$\begin{aligned}&\sum_{i=1}^n \mathbb{E}_1\left|\frac{Z_{i,n}-\mathbb{E}_1(Z_{i,n})}{\sqrt{n\cdot V_1(Z_{1,n})}}\right|^4\\ &\quad\le \frac{n\mathbb{E}_1|Z_{1,n}|^4}{(c'_v)^2n^2\kappa_n^2}\le \frac{n\int[G(\cdot,m)\star f_\varepsilon]^4(x-x_0)\,\mathrm{d}x(1+\mathrm{o}(1))}{(c'_v)^2\ln^2(n)}\\ &\quad\le \frac{nc\int|G^*(u,m)f^*_\varepsilon(u)|^2\,\mathrm{d}u(\int|G^*(u,m)f^*_\varepsilon(u)|\,\mathrm{d}u)^2}{(c'_v)^2\ln^2(n)}\\ &\quad\le c\frac{\ln(n)\cdot m^{-2\underline{b}-2\gamma+1}}{\ln^2(n)}=\mathrm{o}(1),\end{aligned}$$

as $n\to\infty$ and since $\underline{b}>1/2$. Next we apply Lyapunov's central limit theorem for triangular arrays, see Petrov [27], to get $\mathbb{P}_1(U_n\ge u_n)\ge 1-\epsilon$, as, when $n\to+\infty$,

$$0\ge u_n=\frac{\ln(\tau)+\kappa_n}{\sqrt{c_v\kappa_n}}=\frac{-(2\gamma+1)/(2\underline{b}+2\gamma)+c_\chi c_0}{\sqrt{c_vc_\chi c_0}}\sqrt{\ln(n)}\to-\infty.$$

(2) Case $\alpha,r>0$ and $\rho=0$. Without loss of generality we consider $b=0$.

In this case, take some $a\in[\underline{a},\overline{a}]$ and $g_0$ belonging to $\mathcal{S}(a,\overline{r},\overline{L}/2)$ such that $g_0>0$ and $g_0(x)\ge c|x|^{-2}$ as $|x|\to\infty$. Let us consider a function $G$ as for the case 1 such that $G^*$ is three-times continuously differentiable having the property

$$\frac{I(\pi/2 \le |u| \le 3\pi/4)}{c(1+u^4)} \le G^*(u) \le \frac{I(\pi/4 \le |u| \le \pi)}{c(1+u^4)}.$$

Next, $g_{1,n}(x) = g_0(x) + \sqrt{c_0 \ln\ln n/n}\, m^{\gamma+1/2} G(m(x-x_0))$, where $m$ is such that

$$c_0 \frac{\ln\ln n}{n} m^{2\gamma+\underline{r}-1} \exp(2\underline{a}(\pi m)^{\underline{r}}) \le 2\pi\overline{L}/2. \tag{37}$$

Note that this gives a first-order approximation of $m = (\log n/(2\underline{a}))^{1/\underline{r}}$. Then, similarly to the case 1, $g_{1,n}$ is a proper density function as soon as $n$ is large enough and for some $M > 0$ we have $f^Z_{1,n}(x) = g_{1,n} * f_\varepsilon(x) \ge C|x|^{-2}$ for all $|x| \ge M$.

By using (37), we get that $g_{1,n}$ belongs to $S(\underline{a}, \underline{r}, \overline{L})$ for any $a \ge \underline{a}$. Next, $|g_{1,n}(x_0) - g_0(x_0)|/\phi_{n,\underline{a},\underline{r}} = c_0|G(0)| > 0$ and we get, in the same way as for case 1,

$$\begin{aligned}
\chi^2(f^Z_0, f^Z_{1,n}) &= c_0 \frac{\ln\ln n}{n} m^{2\gamma+1} \int \frac{[G(m(\cdot - x_0)) \star f_\varepsilon]^2(x)}{f^Z_{1,n}(x)}\, \mathrm{d}x \\
&\le c_0 \frac{\ln\ln n}{n} m^{2\gamma+1} c_1 \int [G(m(\cdot - x_0)) \star f_\varepsilon]^2(x)\, \mathrm{d}x (1+\mathrm{o}(1)) \\
&\le c_0 c_\chi \frac{\ln\ln n}{n} =: \kappa_n.
\end{aligned}$$

Let us choose $c_0$ small such that $c_0 c_\chi < (\overline{r} - \underline{r})(2\gamma+1)/(\overline{r}\underline{r})$ and let $\xi$ and $\tau$ be defined by

$$c_0 c_\chi < \xi < \frac{\overline{r} - \underline{r}}{\overline{r}\underline{r}}(2\gamma+1) \text{ and } \tau = \ln(n)^{-\xi}.$$

On the one hand, this implies $\tau q_n^2 \to \infty$ with $n$. On the other hand, after checking again that Lyapunov's central limit theorem holds in this case we get

$$\mathbb{P}_1(\mathrm{d}\mathbb{P}_0/\mathrm{d}\mathbb{P}_1 \ge \tau) \ge \mathbb{P}_1(U_n \ge u_n) \ge 1 - \epsilon,$$

as $u_n = (-\ln(\tau) + n\kappa_n)(c_v n \kappa_n)^{-1/2} = (-\xi + c_0 c_\chi)(c_v c_0 c_\chi)^{-1/2}\sqrt{\ln\ln(n)} \to -\infty$.

(3) Case $r > 0$, $0 < \rho \le 1$ and $r \in [\underline{r}, \overline{r}]$ such that $\underline{r} \ge \rho$. Without loss of generality we consider $b = 0$.

As in the second case, take some $a \in [\underline{a}, \overline{a}]$ and $g_0$ belonging to $\mathcal{S}(a, \overline{r}, \overline{L}/2)$ such that $g_0 > 0$ and $g_0(x) \ge c|x|^{-2}$ as $|x| \to \infty$. Let also $G$ be a function such that $G^*$ is three-times continuously differentiable with a bounded first derivative and having the property

$$I(\pi/2 \le |u| \le 3\pi/4) \le G^*(u) \le I(\pi/4 \le |u| \le \pi).$$

Next, define $g_{1,n}$ via its Fourier transform

$$g^*_{1,n}(u) = g^*_0(u) + c_0 \frac{\mathrm{e}^{-\alpha(\pi m)^\rho}}{\sqrt{n}} m^{\rho-1/2} \mathrm{e}^{2\alpha|u|^\rho} G^*(|u|^\rho - (\pi m)^\rho) \mathrm{e}^{\mathrm{i}ux_0},$$

where $m$ is the solution of the equation

$$2\underline{a}(\pi m)^{\underline{r}} + 2\alpha(\pi m)^{\rho} = \log n - (\log\log n)^2. \tag{38}$$

We stress the fact that $m$ is no longer a scaling parameter of the function $G$ in this construction.

Again, as previously, we can check that $g_{1,n}$ is a proper probability density, as soon as $n$ is large enough, and that for some $M > 0$ we have $f_{1,n}^Z(x) \geq C|x|^{-2}$ for all $|x| \geq M$.

Let us check that $g_{1,n}$ belongs to $S(\underline{a}, \underline{r}, \overline{L})$. It is enough to bound from above

$$\begin{aligned}
&(2\pi n)^{-1} \int c_0^2 \mathrm{e}^{-2\alpha(\pi m)^{\rho}} m^{2\rho-1} \mathrm{e}^{4\alpha|u|^{\rho}} |G^*(|u|^{\rho} - (\pi m)^{\rho})|^2 \mathrm{e}^{2\underline{a}|u|^{\underline{r}}} \,\mathrm{d}u \\
&\quad \leq (2\pi n)^{-1} c_0^2 m^{2\rho-1} \mathrm{e}^{-2\alpha(\pi m)^{\rho}} \int_{\pi/4 \leq |u|^{\rho} - (\pi m)^{\rho} \leq 3\pi/4} \mathrm{e}^{4\alpha|u|^{\rho} + 2\underline{a}|u|^{\underline{r}}} \,\mathrm{d}u \\
&\quad \leq (2\pi n)^{-1} c_0^2 c_1 m^{2\rho-1} \mathrm{e}^{-2\alpha(\pi m)^{\rho}} (\pi m)^{1-\underline{r}} \mathrm{e}^{4\alpha(\pi m)^{\rho} + 2\underline{a}(\pi m)^{\underline{r}}} \\
&\quad \leq c_0^2 c_2 n^{-1} m^{2\rho - \underline{r}} \mathrm{e}^{2\underline{a}(\pi m)^{\underline{r}} + 2\alpha(\pi m)^{\rho}},
\end{aligned}$$

which tends to 0 when $m$ is defined by (38). Next,

$$\begin{aligned}
|g_{1,n}(x_0) - g_0(x_0)| &= (2\pi\sqrt{n})^{-1} \left| \int c_0 \mathrm{e}^{-\alpha(\pi m)^{\rho}} m^{\rho-1/2} \mathrm{e}^{2\alpha|u|^{\rho}} G^*(|u|^{\rho} - (\pi m)^{\rho}) \,\mathrm{d}u \right| \\
&\geq c_0 m^{\rho-1/2} \frac{\mathrm{e}^{-\alpha(\pi m)^{\rho}}}{2\pi\sqrt{n}} \int_{\pi/2 \leq |u|^{\rho} - (\pi m)^{\rho} \leq \pi} \mathrm{e}^{2\alpha|u|^{\rho}} \,\mathrm{d}u \\
&\geq c_0 c_3 m^{1/2} \frac{\mathrm{e}^{\alpha(\pi m)^{\rho}}}{2\pi\sqrt{n}}
\end{aligned}$$

and we can check similarly to Butucea and Tsybakov [6] that for $m$ solution of (38) this sequence is equivalent to $\phi_{n,\underline{a},\underline{r}}$ when $n \to \infty$. Finally

$$\begin{aligned}
\chi^2(f_0^Z, f_{1,n}^Z) &= c_0^2 \int [(g_{1,n} - g_0) \star f_{\varepsilon}]^2(x) / f_{1,n}^Z(x) \,\mathrm{d}x \\
&\leq c_0^2 \left\{ \int_{|x| \leq M} [(g_{1,n} - g_0) \star f_{\varepsilon}]^2(x) \,\mathrm{d}x + \int_{|x| > M} x^2 [(g_{1,n} - g_0) \star f_{\varepsilon}]^2(x) \,\mathrm{d}x \right\},
\end{aligned}$$

say $T_1 + T_2$. Then

$$\begin{aligned}
T_1 &\leq c_0^2 c_4 n^{-1} \mathrm{e}^{-2\alpha(\pi m)^{\rho}} m^{2\rho-1} \int |G^*(|u|^{\rho} - (\pi m)^{\rho}) f_{\varepsilon}^*(u)|^2 \,\mathrm{d}u \\
&\leq c_0^2 c_5 n^{-1} \mathrm{e}^{-2\alpha(\pi m)^{\rho}} m^{2\rho-1} \int_{\pi/4 \leq |u|^{\rho} - (\pi m)^{\rho} \leq 3\pi/4} \mathrm{e}^{2\alpha|u|^{\rho}} \,\mathrm{d}u = c_0^2 c_6 (\pi m)^{\rho}/n.
\end{aligned}$$

Moreover, under the additional assumption (22) that $|\partial f_\varepsilon^*(u)/\partial u| \le \mathrm{O}(1)|u|^{\rho-1}\exp(-\alpha|u|^\rho)$ as $|u|\to\infty$,

$$T_2 \le c_0^2 c_7 n^{-1}\mathrm{e}^{-2\alpha(\pi m)^\rho} m^{2\rho-1}\int\left|\frac{\partial}{\partial u}[G^*(|u|^\rho-(\pi m)^\rho)f_\varepsilon^*(u)]\right|^2\mathrm{d}u$$

$$\le c_8 n^{-1}\mathrm{e}^{-2\alpha(\pi m)^\rho} m^{2\rho-1}\int_{\pi/4\le|u|^\rho-(\pi m)^\rho\le 3\pi/4}|u|^{2(\rho-1)}\mathrm{e}^{2\alpha|u|^\rho}\,\mathrm{d}u \le c_9\frac{(\pi m)^{3\rho-2}}{n}=\mathrm{o}(T_1),$$

for $\rho\le 1$ and $n$ large enough. Thus $\chi^2(f_0^Z, f_{1,n}^Z)\le c_0^2 c_\chi(\pi m)^\rho/n =: \kappa_n$.

Let $c_0$ be small such that $c_0^2 c_\chi < 2\alpha$ and let $\xi$ and $\tau$ be defined by $c_0^2 c_\chi<\xi<2\alpha$ and $\tau=\mathrm{e}^{-\xi(\pi\ln(n)/(2\underline{a}))^{\rho/\underline{r}}}$. We have $\tau\phi^2_{n,\underline{a},\underline{r}}/\phi^2_{n,\overline{a},\overline{r}}\ge(\ln(n))^A\exp((-\xi+2\alpha)(\ln(n)/(2\underline{a}))^{\rho/\underline{r}}+B(\ln(n))^C)$ tends to infinity for some real numbers $A$, $B$, $C$, as $C<\rho/\underline{r}$ and $\xi<2\alpha$. We check that Lyapunov's theorem holds and that

$$u_n=\frac{-\ln(\tau)+n\kappa_n}{\sqrt{c_v n\kappa_n}}=\frac{-\xi(\pi\ln(n)/(2\underline{a}))^{\rho/\underline{r}}+c_0^2c_\chi(\pi m)^\rho}{c_0\sqrt{c_v c_\chi}(\pi m)^{\rho/2}}\to-\infty$$

with $n$, as $m$ defined by (38) is larger than $(\ln(n)/(2\underline{a}))^{1/\underline{r}}$.

The proof that $\varphi_n$ is the minimax rate of estimation in this case repeats the proof of (3) with modified choice of $g_{1,n}$ via its Fourier transform

$$g_{1,n}^*(u)=g_0^*(u)+c_0\frac{\mathrm{e}^{-\alpha(\pi m)^\rho}}{\sqrt{n}}m^{(\rho-1)/2}\mathrm{e}^{2\alpha|u|^\rho}G^*(|u|^\rho-(\pi m)^\rho)\mathrm{e}^{\mathrm{i}ux_0},$$

where $m$ is the solution of equation (38).

This gives the rate $|g_{1,n}(x_0)-g_0(x_0)|\ge c_0c_3m^{-(\rho-1)/2}\mathrm{e}^{\alpha(\pi m)^r}/\sqrt{n}$, which is equivalent to $V_{\breve{m}}$ for $n$ large enough and $n\chi^2(f_0^Z, f_{1,n}^Z)\le c_0^2c_6+c_9m^{2\rho-2}\le c_0^2c_\chi$. Thus, the rate $\varphi_n$ is a minimax rate of convergence for $r\ge\rho$, $\rho\le 1$.